\documentstyle[leqno,12pt]{article}
\newlength{\abstractwidth}
\renewcommand{\thefootnote}{\fnsymbol{footnote}}
\renewcommand{\thanks}[1]{\footnote{#1}} 
\newcommand{\starttext}{ \setcounter{footnote}{0}
\renewcommand{\thefootnote}{\arabic{footnote}}}

\newcommand{\be}{\begin{equation}}
\newcommand{\bea}{\begin{eqnarray}}
\newcommand{\eea}{\end{eqnarray}} \newcommand{\ee}{\end{equation}}
 
 \def\ba{\begin{eqnarray}}
\def\ea{\end{eqnarray}}

\def\o{\omega}

\def\tr{{\rm tr}}

\def\log{\,{\rm log}\,}
\def\exp{\,{\rm exp}\,}

\def\o{\omega}

\def\o{\omega}

\def\ge{\geq}
\def\le{\leq}

\def\p{\partial}

\newcommand{\bk}[1]{\Big(#1\Big)}

\def\[{{\bf [}}
\def\]{{\bf ]}}

\def\ddbar{i\p\bar\p}

\def\ric{{\rm Ric}}
\def\mathbb{\bf}

\def\vol{\mathrm{Vol}}

\newcommand{\neweqref}[1]{(\ref{#1})}

\begin{document}
\starttext \baselineskip=18pt \setcounter{footnote}{0}
\newtheorem{theorem}{Theorem}
\newtheorem{lemma}{Lemma}
\newtheorem{corollary}{Corollary}
\newtheorem{definition}{Definition}
\newtheorem{conjecture}{Conjecture}
\newtheorem{proposition}{Proposition}
\newcommand{\beqref}[1]{(\ref{#1})}


\begin{center}
{\Large \bf Green's functions and complex Monge-Amp\`ere equations
\footnote{Work supported in part by the National Science Foundation under grant DMS-1855947.}}

\medskip
\centerline{Bin Guo, Duong H. Phong, and Jacob Sturm}

\medskip

\begin{abstract}

{\footnotesize Uniform $L^1$ and lower bounds are obtained for the Green's function on compact K\"ahler manifolds. Unlike 
 in the classic theorem of Cheng-Li for Riemannian manifolds, the lower bounds do not depend directly on the Ricci curvature, but only on integral bounds for the volume form and certain of its derivatives. In particular, a uniform lower bound for the Green's function on K\"ahler manifolds is obtained which depends only on a lower bound for the {\it scalar} curvature and on an $L^q$ norm for the volume form for some $q>1$. The proof relies on auxiliary Monge-Amp\`ere equations, and is fundamentally non-linear. The lower bounds for the Green's function imply in turn $C^1$ and $C^2$ estimates for complex Monge-Amp\`ere equations with a sharper dependence on the function on the right hand side.}

\end{abstract}

\end{center}

\baselineskip=15pt
\setcounter{equation}{0}
\setcounter{footnote}{0}

\section{Introduction}
\setcounter{equation}{0}

A fundamental theorem in Riemannian geometry is the inequality of Cheng-Li \cite{CL}, which provides conditions for a uniform lower bound for the Green's function. More precisely, let $(X,g)$ be a compact Riemannian manifold, and define the Green's function $G(x,y)$ as the unique function $X\times X\to {\bf R}\cup\{\infty\}$ satisfying the conditions 
\begin{equation}\label{eqn:Green}\Delta_{g, y} G(x, y) = -\delta_x(y) + \frac{1}{{\mathrm{Vol}}_g(X)},\quad \int_X G(x,y) dV_g(y) = 0,\end{equation}
where $\delta_x(y)$ denotes the Dirac function at the point $x$. It is well-known that $G(x,y) = G(y,x)$ and $\Delta_{g,x} G(x,y) = \Delta_{g, y} G(x,y)$ for any $x\neq y$ . 
Assume that the Ricci curvature of $(X,g)$ satisfies $\ric(g)\geq -\kappa$ for some constant $\kappa$. Then Cheng-Li \cite{CL} prove that there is a constant $C>0$ depending only on the dimension of $X$ and $\kappa$ such that
\begin{equation}\label{eqn:Cheng Li}G(x,y) \ge - C \frac{{\mathrm{diam}_g}(X)^2}{{\vol}_g(X) }, \quad \forall x, y\in X.\end{equation}
Here ${\mathrm{diam}}_g(X)$ is the diameter of $(X,g)$ and ${\vol}_g(X)$ is its volume. The dependence of this inequality on a lower bound $\kappa$ for the Ricci curvature is crucial, and it does not seem possible in Riemannian geometry to lessen this dependence in any significant way.

\medskip
A first primary goal of the present paper is to show that, in the K\"ahler setting, lower bounds for the Green's function can actually be established without assumptions on lower bounds for the Ricci curvature. Rather, we assume integral bounds on the volume form and some specific derivatives. As we shall explain below, for our purposes, integral bounds are far superior to pointwise bounds. It turns out that several distinct sets of assumptions can guarantee lower bounds for the Green's function, and we shall describe them in detail later in \S 2 (Theorem 2). While these sets of assumptions may be difficult to assess at a glance, an easy comparison with the theorem of Cheng-Li can be obtained by observing that Theorem 2 implies in particular that the Green's function is bounded from below if the scalar curvature is bounded from below, and the $L^q$ norm of the volume form is bounded for some $q>1$
(see Corollary \ref{cor CL}).
It may also be worth stressing that, while our goal of establishing lower bounds for the Green's function is the same as Cheng-Li's, our method of proof is completely different. It builds repeatedly on the key idea in \cite{GPT} of comparison with an auxiliary Monge-Amp\`ere equation, and even though the Green's function is the solution of a linear partial differential equation, our method is fundamentally non-linear.

\medskip
We come now to the second primary goal of this paper, which is also a major motivation for the first, and which is sharp a priori estimates for general complex Monge-Amp\`ere equations.
The first estimates were obtained in 1976 by S.T. Yau in his seminal paper \cite{Y} solving the Calabi conjecture. However, a new generation of problems in complex geometry is leading to more complicated Monge-Amp\`ere equations, which can be degenerate or singular in many different senses. Thus ever sharper versions of a priori estimates are needed, as each improvement usually has significant geometric consequences.

\medskip
The first sharp form of $C^0$ estimates for the complex Monge-Amp\`ere equation was obtained by Kolodziej \cite{K}, using pluripotential theory. Kolodziej's estimates were extended to the important case of degenerating background K\"ahler metrics by Eyssidieux, Guedj, and Zeriahi \cite{EGZ} and Demailly and Pali \cite {DP}.  More general $C^0$ estimates using the theory of envelopes can be found in \cite{GL, GL1}. Another approach to $C^0$ estimates, using PDE methods, was introduced very recently in \cite{GPT}. This method can also apply to nef classes \cite{GPTW}, and lead to many sharp estimates, including stability estimates \cite{GPTa}, diameter estimates \cite{GPTW1}, and non-collapse estimates \cite{GS}. We have just seen it applied to lower bounds for Green's functions in the first part of this paper. On the other hand, while $C^1$ and $C^2$ estimates were extended to the case allowing a divisor, building on ideas of Tsuji \cite{Ts}, Blocki \cite{B}, and \cite{PS, PSoS}, they still require very restrictive conditions on the right hand side of the equation, such as pointwise lower or upper bounds, and in the case of $C^2$ estimates, also a bound on its Laplacian. A lower or upper bound assumption is a particularly severe constraint, as it may rule out equations which are degenerate or have singularities.

\medskip
Thus our second primary goal in this paper is to establish sharp $C^1$ and $C^2$ estimates for the complex Monge-Amp\`ere equation which depend essentially only on integral bounds for the right hand side. In order to do so, we cannot apply the standard maximum principle to the elliptic differential inequalities satisfied by the derivatives of the solution of the Monge-Amp\`ere equation. Rather, we apply instead the new lower bounds for the Green's function obtained in the first part of the paper.

\medskip
We now state precisely our main results. In view of many applications in complex geometry, it is important to obtain estimates which remain uniform as the K\"ahler class may degenerate, so we consider the following set-up, which includes both cases of fixed and degenerating K\"ahler classes as special cases.
Let $(X,\omega_X)$ be a compact K\"ahler manifold with dimension $n$. Suppose $\chi$ is a $d$-closed $(1,1)$-form on $X$ such that its cohomology class $[\chi]\in H^{1,1}(X,{\mathbb C})$ is {\em nef} and big, i.e. $[\chi]$ lies in the closure of the K\"ahler cone of $X$ and $\int_X \chi^n>0$. So for each $t>0$, $[\chi + t\omega_X]$ is a K\"ahler class. For any fixed $t\in (0,1]$ and any K\"ahler metric $\omega_t\in [\chi + t \omega_X]$, we define the function $F_{\omega_t}$ to be the {\em logarithmic of the relative volume form of $\omega_t$} with respect to the fixed volume form $\omega_X^n$, that is
\begin{equation}\label{eqn:relative volume}
F_{\omega_t} = \log \Big(\frac{\omega_t^n/V_t}{\omega_X^n/V} \Big)
\end{equation}
where $V_t = \int_X \omega_t^n = \int_X(\chi + t\omega_X)^n>0$ is the volume of the K\"ahler class $[\chi + t\omega_X]$, and $V = \int_X \omega_X^n$ is the volume of the fixed metric $\omega_X$. Note that $V_0 = \int_X \chi^n>0$ under our assumptions.

\smallskip

Fix $p>n$. We denote the $p$-th entropy of the K\"ahler metric $\omega_t$ by
$${\mathrm {Ent}}_p(\omega_t) = \frac{1}{V_t} \int_X |F_{\omega_t}|^p {\omega_t^n} = \frac{1}{V} \int_X |F_{\omega_t}|^p e^{F_{\omega_t}} \omega_X^n.$$
For any $N>0$, we define a subset of the space of K\"ahler metrics in $[\chi + t\omega_X]$ by
$${\mathcal M}_t(N,p) = \{\omega_t\in [\chi + t\omega_X]|~ {\mathrm {Ent}}_p(\omega_t) \le N \}.$$

Our first main theorem shows that the lower bound of $G_t$ is equivalent to the $L^1(X,\omega_t^n)$-norm of $G_t$ for $\omega_t\in {\mathcal M}_t(N,p)$.

\begin{theorem}\label{thm:main1}
Given $p>n$ and $N>0$, there is a constant  
%
$C>0$ depending only on $n, p, \chi, \omega_X$ and $N$,
such that for any $t\in (0,1]$ and any K\"ahler metric $\omega_t\in {\mathcal M}_{t}(N, p)$, the Green's function $G_t$ associated to $\omega_t$ satisfies
$$\frac{1}{2 V_t} \| G_t(x,\cdot)\|_{L^1(X,\omega_t^n)} \le - \inf_X G_t(x,\cdot)\le C(1 + \| G_t(x,\cdot)\|_{L^1(X,\omega_t^n)}),$$ for any $x\in X$.
\end{theorem}

\medskip

We remark that $G_t(x,y)$ satisfies the  asymptotic behavior \cite{SY}
 \begin{equation}\label{eqn:G asymp}G_t(x,y)\sim d_{\omega_t}(x,y)^{-2n + 2}, \mbox{ if }n\ge 2\end{equation} 
 and $G_t(x,y)\sim - \log d_{\omega_t}(x,y)$ if $n =1$, when $x$ is close to $y$. Here $d_{\omega_t}(x,y)$ denotes the geodesic distance of $x$ and $y$ under $\omega_t$. Thus an upper bound for $G_t(x,y)$ cannot be expected to hold. 

\bigskip

For a smaller class of K\"ahler metrics than ${\mathcal M}_t(N,p)$, we will show that the $L^1$-norms of the Green's function are uniformly bounded, hence by Theorem \ref{thm:main1}, we also have a pointwise lower bound on the Green's functions.

\smallskip

Henceforth we require that $\chi$ is {\em nonnegative} and $V_0 = \int_X\chi^n>0$. Let $F_{\omega_t}$ be associated with a K\"ahler metric $\omega_t\in [\chi + t\omega_X]$ as in \neweqref{eqn:relative volume}. For $\epsilon>0$, $N>0$ and $\gamma\ge 1$, we denote for each $t\in (0,1]$
\begin{equation}\label{eqn:class 1}{\mathcal M}'_t(N,\epsilon,\gamma) = \Big\{\omega_t\in [\chi + t\omega_X]|~ \frac{1}{V} \int_X e^{(1+\epsilon) F_{\omega_t}} \omega_X^n \le N  \mbox{ and }\sup_X e^{-F_{\omega_t}}\le \gamma  \Big \}.\end{equation}
A second class of metrics for $t\in (0,1]$ is given by
\bea\label{eqn:class 2}{\mathcal M}_t''(N,\epsilon,\gamma) & = &\Big\{\omega_t\in [\chi + t\omega_X]|~ \frac{1}{V} \int_X e^{(1+\epsilon) F_{\omega_t}} \omega_X^n \le N \\
& &\qquad  \mbox{ and }\nonumber \int_X (e^{-F_{\omega_t}} + |\Delta_{\omega_X} e^{-F_{\omega_t}}|) \omega_X^n\le  \gamma \Big \}.
\eea
The case of K\"ahler metrics in a fixed K\"ahler class $\omega\in [\omega_X]$ can be considered as a special case of the above more general set-up by taking
$\chi = \omega_X/2$ and $t = 1/2$. We can also consider the class: 
\bea\label{eqn:class 3}\tilde{\mathcal M}''(N,\epsilon,\gamma) & = &\Big\{\omega \in [\omega_X]|~ \frac{1}{V} \int_X e^{(1+\epsilon) F_{\omega}} \omega_X^n \le N \\
& &\qquad  \mbox{ and }\nonumber \int_X (e^{-F_{\omega}} + |\nabla _{\omega_X} e^{-F_{\omega}}|_{\omega_X}^2) \omega_X^n\le  \gamma \Big \}.
\eea
Abusing notations, when we  write $\omega_t\in \tilde{\mathcal M}''(N,\epsilon,\gamma)$, we mean that $t = 1/2$ and $\chi = \omega_X/2$, which corresponds to the case of complex Monge-Amp\`ere equations with the fixed background metric $\omega_X$.

\smallskip

It follows easily from calculus inequalities that the sets of metrics in \neweqref{eqn:class 1}, \neweqref{eqn:class 2}, \neweqref{eqn:class 3}  are contained in ${\mathcal M}_t(\tilde N, p)$ for suitable $\tilde N> N$ and $p>n$. Hence Theorem \ref{thm:main1} applies to the metrics in these sets. Our second theorem is:

\medskip
\begin{theorem}\label{thm:main2}
Given any $\epsilon>0$, $N>0$ and $\gamma\in (0,1)$, for each $t\in (0,1]$ and each K\"ahler  metric $\omega_t\in {\mathcal M}_t'(N,\epsilon,\gamma) \cup  {\mathcal M}_t''(N,\epsilon,\gamma) \cup \tilde {\mathcal M}''(N,\epsilon,\gamma)$, let $G_t$ be the Green's function associated with $(X, \omega_t)$. Then the following hold:

\smallskip

 (i) There is a constant $C>0$ which depends on $n, \epsilon, \omega_X,\chi$ and $N, \gamma$ such that for any $x\in X$
\begin{equation}\label{eqn:Green bound}
\| G_t(x,\cdot)\|_{L^1(X,\omega_t^n)} = \int_X |G_t(x,\cdot)|\omega_t^n \le C,
\end{equation}
and  \begin{equation}\label{eqn:Green lower bound}\inf_{y\in X} G_t(x,y)\ge - C.\end{equation}

\smallskip

 (ii) For any given  $\delta\in (0, \frac{2n}{2n-1})$, there is a constant $C_\delta>0$ depending additionally on $\delta$ such that  $G_t$ satisfies
$$ \int_X |G_t(x,\cdot) |^{\frac{n}{n-1} - \delta} \omega_t^n + \int_X |\nabla G_t(x,\cdot)|_{_{\omega_t}}^{\frac{2n}{2n-1} - \delta} \omega_t^n \le C_\delta, $$ 
for any fixed $x\in X$.

\smallskip

Moreover, when $n = 1$, (i) and (ii) hold for any $\omega_t$ with $\| e^{F_{\omega_t}}\|_{L^{1+\epsilon}}\le N$, and no extra conditions on $e^{-F_{\omega_t}}$ as in \neweqref{eqn:class 1}, \neweqref{eqn:class 2} and \neweqref{eqn:class 3} are needed.
\end{theorem}


We now turn to the application of the estimates of Green's functions in Theorem  \ref{thm:main2}. Again we assume $\chi$ is {\em nonnegative} and the class $[\chi]$ is big. Then by Kodaira's lemma, there is an effective divisor $D$ on $X$ such that 
$$\chi - \varepsilon_0 \ric(h_D) \ge \delta_0 \omega_X$$ for suitable positive constants $\varepsilon_0$ and $\delta_0$ which we will fix throughout the paper, where $h_D$ is a Hermitian metric on the line bundle $[D]$ associated with $D$. Let $s_D\in {\mathcal O}_X(D)$ be a holomorphic section defining $D$ such that
$$
{\rm sup}_X|s_D|_{h_D}^2=1.
$$
Let $\varphi_t$ be the K\"ahler potential of the K\"ahler metric $\omega_t \in [\chi + t\omega_X]$, i.e. $\omega_t = \chi + t\omega_X + \ddbar \varphi_t$. We first derive an estimate on the gradient of $\varphi_t$ with respect to the fixed metric $\omega_X$, for $\omega_t$ in the sets of K\"ahler metrics in Theorem \ref{thm:main2}.
%
\smallskip

\begin{theorem}\label{thm:gradient}
Given $N, \epsilon,\gamma\in (0,1)$, and $p>n$, for each  $t\in (0,1]$ and any $\omega_t\in {\mathcal M}_t'(N,\epsilon,\gamma) \cup  {\mathcal M}_t''(N,\epsilon,\gamma) \cup \tilde {\mathcal M}''(N,\epsilon,\gamma)$, the following estimate on $|\nabla \varphi_t|_{\omega_X}^2$ holds. There is a constant $C>0$ depending on $n, \epsilon, \chi,\omega_X, N,\gamma$, $p$ and  $\int_X |\nabla F_{\omega_t}|_{\omega_X}^p e^{F_{\omega_t}} \omega_X^n$ such that away from $D$
$$ |\nabla \varphi_t|_{\omega_X}^2 \le \frac{C}{|s_D|_{h_D}^{2 A}},$$
where 
$A>0$ depends only on $n, \epsilon, \chi, \omega_X, N$ and $\gamma$.
\end{theorem}%

\smallskip

We observe that gradient estimates for solutions to certain complex Monge-Amp\`ere equations had been obtained in \cite{B, PS}, 
 but they required pointwise bounds on $|\nabla F_{\omega_t}|_{\omega_X}$. In \cite{CH, GPT1}, the gradient is shown to depend on certain $L^p$ bound of $|\nabla F_{\omega_t}|_{\omega_X}$ for $p\ge 2n$. 
 Theorem \ref{thm:gradient} implies that the gradient estimate depends only on integral bound of $e^{F_{\omega_t}}$ and $L^p$-bound of $|\nabla F_{\omega_t}|_{\omega_X}$ for smaller $p$ which turns out to be sharp. In particular, the gradient estimate continues to hold in some situations even when $e^{F_{\omega_t}}$ has zeros or is unbounded. Theorem \ref{thm:gradient} also provides immediately a gradient estimate for solutions to complex Monge-Amp\`ere equations with a fixed background metric (i.e. when $\omega\in  \tilde {\mathcal M}''(N,\epsilon,\gamma)$, cf. Corollary \ref{corollary for gradient} in \S \ref{section 4} below), in which case we can take $s_D\equiv 1$, $h_D \equiv 1$ and $D$ to be trivial. {Example 3.1} shows that the gradient estimate may fail if $p<n$, so the assumption on $p>n$ in Theorem \ref{thm:gradient} is sharp.

\smallskip

With regard to the second-order derivatives, we have the following counterpart of Yau's $C^2$ estimate \cite{Y}.

\begin{theorem}\label{thm:C2}
Under the same setup as in Theorem \ref{thm:gradient} with $p>2n$, 
there is a constant $C>0$ depending on $n, \epsilon, \chi,\omega_X, N,\gamma$, $p$ and  $\int_X |\nabla F_{\omega_t}|_{\omega_X}^p e^{F_{\omega_t}} \omega_X^n$ such that away from $D$ the following holds

$$ |\ddbar  \varphi_t|_{\omega_X}^2 \le \frac{C}{|s_D|_{h_D}^{2 B}},$$
where 
$B>0$ depends only on $n, \epsilon, \chi, \omega_X,N$ and $\gamma$.
\end{theorem}
We stress that the above $C^2$-estimate of $\varphi_t$ is independent of the pointwise estimates on the second order derivatives of $F_{\omega_t}$, and it only depends on certain $L^p$-bound of $|\nabla F_{\omega_t}|_{\omega_X}$. Theorem \ref{thm:C2} improves in particular on the known estimates \cite{Y} for complex Monge-Amp\`ere equations with a fixed background metric (cf. Corollary \ref{cor 3.2}). We remark that the exponent $p>2n$ in Theorem \ref{thm:C2} is also sharp, as  {Example 3.2} shows that the estimates may fail if $p< 2n$. 

\smallskip

By utilizing the estimates of the Green's functions, we can also improve the $C^3$ estimates of complex Monge-Amp\`ere equations (cf. Theorem \ref{thm:C3}), which generalize the known ones in \cite{Y, PSS} by a weaker dependence of the function $e^F$ on the right-hand side.


\section{Proof of Theorem \ref{thm:main1}} \setcounter{equation}{0}

Given any $t\in (0,1]$, we fix an arbitrary K\"ahler metric $\omega_t\in {\mathcal M}_{t}(N, p)$.  It is clear from the $\ddbar $-lemma that \neweqref{eqn:relative volume} is equivalent to the following 
complex Monge-Amp\`ere equation with $\omega_t = \chi +  t\omega_X + \ddbar \varphi_t$
\begin{equation}\label{eqn:MA}
(\chi + t \omega_X + \ddbar \varphi_t)^n = c_t e^{F_{\omega_t}} \omega_X^n,\quad \sup_X\varphi_t = 0,
\end{equation}
where we have normalized $\varphi_t$ to make it unique, and we denote $c_t = V_t/V$. The case of a fixed K\"ahler class $\o_X$ corresponds for example to the special case $\chi={1\over 2}\o_X$ and $t={1\over 2}$. For simplicity of notations, we will write $F_{\omega_t}$ as $F$.

\smallskip

Since $\hat \omega_t: = \chi+ t\omega_X$ may not necessarily be positive, we introduce the following $\hat\omega_t$-plurisubharmonic (PSH) function.
\begin{definition}
For any $t\in (0,1]$, we denote the {\em envelope} associated to the $(1,1)$-form $\hat \omega_t$
$${\mathcal {V}}_t = \sup\{v\in PSH(X,\hat\omega_t)| ~ v\le 0\}.$$
\end{definition}%
Note that if $\chi\ge 0$ is a nonnegative $(1,1)$-form, ${\mathcal V}_t \equiv 0$ for any $t\in (0,1]$.

\smallskip

The following $L^\infty$ estimate for the family of solutions $\varphi_t$ to the equation \neweqref{eqn:MA} is proved in \cite{GPTW} (see also \cite{K, FGS}). 
\begin{lemma}[\cite{GPTW}]\label{lemma GPTW}
There is a uniform constant ${C}_0={C}_0(n, p, \chi, \omega_X, {\mathrm{Ent}}_p(\omega_t))>0$ such that 
$$\sup_X | \varphi_t - {\mathcal V}_t|\le {C_0},\quad \forall t\in (0,1].$$
\end{lemma}
To prove Theorem \ref{thm:main1}, we need the following mean-value type inequality for functions satisfying certain linear differential inequality.
\begin{lemma}\label{lemma key}
Suppose $v\in L^1(X,\omega_t^n)$ is a function that satisfies $\int_X v \omega_t^n= 0$ and 
\begin{equation}\label{eqn:assumption 1}
v\in C^2(\Omega_{-1}),\quad \Delta_{\omega_t} v \ge -a \mbox{ in }\Omega_{0}
\end{equation}
for some $a> 0$ and $\Omega_s = \{v> s\}$ is the super-level set of $v$. Then there is a constant $C>0$ depending only on $n, p, \chi, \omega_X, {\mathrm{Ent}}_p(\omega_t)$ and $a$ such that 
$$\sup_X v\le C(1+ \| v\|_{L^1(X,\omega_t^n)}).$$
\end{lemma}
We emphasize that the constant $C>0$ in the lemma above is independent of $t\in (0,1]$. The proof of Lemma  \ref{lemma key} uses the techniques similar to that of the $L^\infty$-estimate for fully nonlinear partial differential equations in \cite{GPT}. The key idea to introduce an auxiliary equation to compare with $v$. The lemma is trivial if $v\equiv 0$, so we assume $v\not\equiv 0$.

\medskip

\noindent{\em Proof. }
We break the proof into four steps. Since the proof is uniform in $t\in (0,1]$, we fix a $t\in (0,1]$. We may assume $\| v\|_{L^1(X,\omega_t^n)}\le V_0$, otherwise, replace $v$ by $\hat v: =V_0 \cdot v /\| v\|_{L^1(X,\omega_t^n)}$ which still satisfies \neweqref{eqn:assumption 1} with the same $a>0$. It suffices to show $\sup_X v\le C$ for some $C>0$ with the dependence as stated in the lemma.

\medskip

\noindent{\bf Step 1}. We fix a sequence of positive smooth functions $\eta_k: {\mathbb R}\to {\mathbb R}_+= (0,\infty)$ such that $\eta_k(x)$ converges uniformly and monotonically decreasingly to the function $x\cdot \chi_{\mathbb R_+}(x)$, as $k\to \infty$. We may choose $\eta_k(x) \equiv 1/k$ for any $x\le -1/2$. For $s\ge 0$ and large $k$, we  consider the following auxiliary complex Monge-Amp\`ere equations
\begin{equation}\label{eqn:aMA}
(\hat \omega_t + \ddbar \psi_{t, k}) ^n = c_t \frac{\eta_k(v - s)}{A_{s, k}} e^F \omega_X^n,\quad \sup_X \psi_{t,k } = 0,
\end{equation}
where 
\begin{equation}\label{eqn:key 1}
A_{s, k} = \frac{c_t}{V_t} \int_X {\eta_k(v - s)} e^F \omega_X^n\to \frac{1}{V} \int_{\Omega_s} (v- s) e^F \omega_X^n=: A_s \mbox{ as }k\to \infty.
\end{equation}
We remark that the right-hand side of \neweqref{eqn:aMA} is smooth and positive, and $[\hat \omega_t]$ is a K\"ahler class, so by Yau's theorem \cite{Y} this equation admits a unique smooth solution $\psi_{t,k}$. We have also assumed that the open set $\Omega_s\neq \emptyset$ so $A_s>0$, otherwise the lemma is already proved. The assumption that $\| v\|_{L^1(X,\omega_t^n)}\le V_0$ implies that $A_s\le 1$, hence $A_{s,k}\le 2$ for large $k$.

\smallskip

\noindent{\bf Step 2.} We denote $\Lambda = C_0 + 1$ where $C_0$ is the constant in Lemma \ref{lemma GPTW}. Consider the function 
$$\Phi: = - \varepsilon (-\psi_{t,k} + \varphi_t + \Lambda)^{\frac n{n+1}} + (v-s),$$ where $\varepsilon>0$ is chosen such that
\begin{equation}\label{eqn:epsilon}
\varepsilon^{n+1} = \Big(\frac{n+1}{n^2} \Big)^n (a+\varepsilon n)^n A_{s,k}.
\end{equation}
It follows easily from $A_{s,k}\le 2$ and equation \neweqref{eqn:epsilon} that \begin{equation}\label{eqn:new epsilon}\varepsilon\le C(n,a) A_{s,k}^{1/(n+1)},\end{equation} for some $C(n,a)>0$ depending only on $n \mbox{ and } a$.
$\Phi$ is a  $C^2$ function on $\Omega_{-1}$ since $v-s$ is so and 
\begin{equation}\label{eqn:key 2}
-\psi_{t,k} + \varphi_t + \Lambda = ({\mathcal V}_t - \psi_{t,k}) + (\varphi_t - {\mathcal V}_t + C_0) + 1\ge 1.
\end{equation}
We claim that $\Phi\le 0$ on $X$. Observe that by the definition of $\Omega_s$ it is clear that $\Phi|_{X\backslash \Omega_s} <0$, so if $\max_{\Omega_s} \Phi \le \sup_{X\backslash \Omega_s} \Phi< 0$, we are done. Otherwise, $\max_{\Omega_s} \Phi > \sup_{X\backslash \Omega_s} \Phi$ and $\Phi$ achieves its maximum at some point $x_0\in \Omega_s$. By maximum principle, $\Delta_{\omega_t} \Phi(x_0)\le 0$. Therefore, we calculate (below we denote $\omega_{t,\psi_{t, k}} = \hat\omega_t + \ddbar \psi_{t,k}$)
\bea
0 &\ge & \Delta_{\omega_t} \Phi(x_0)\nonumber\\
&\nonumber\ge & \frac{\varepsilon n}{n+1} (-\psi_{t,k} + \varphi_t  + \Lambda) ^{-\frac {1}{n+1}} ( \tr_{\omega_t} \omega_{t, \psi_{t, k}} - \tr_{\omega_t} \omega_t ) + \Delta_{\omega_t} v\\
&\ge &\nonumber  \frac{\varepsilon n^2 }{n+1} (-\psi_{t,k} + \varphi_t  + \Lambda) ^{-\frac {1}{n+1}} \Big(\frac{\omega^n _{t,\psi_{t, k}}}{\omega_t^n} \Big)^{1/n} - \frac{\varepsilon n^2}{n+1} - a\\
&\ge &\nonumber  \frac{\varepsilon n^2 }{n+1} (-\psi_{t,k} + \varphi_t  + \Lambda) ^{-\frac {1}{n+1}} \frac{(v-s)^{1/n}}{ A_{s,k} ^{1/n}} - a - \varepsilon n,
\eea
where in the third line we applied the arithmetic-geometric inequality and in the last line we use the equation \neweqref{eqn:aMA}. By the choice of $\varepsilon$ in \neweqref{eqn:epsilon}, it follows by a straightforward calculation that $\Phi(x_0)\le 0$, as claimed. 

\smallskip

\noindent{\bf Step 3.} From $\Phi\le 0$ and \neweqref{eqn:new epsilon} in the previous step , we have $(v-s) A_{s,k}^{-1/(n+1)} \le C_1 (-\psi_{t,k } + \varphi_t + \Lambda)^{n/(n+1)}$ on $X$, for some $C_1>0$ depending only on $n$ and $a$. In particular on $\Omega_s=\{v-s>0\}$ we have by taking $(n+1)/n$-th power
$$\frac{(v-s)^{(n+1)/n}}{A_{s,k}^{1/n}}\le C_1^{(n+1)/n} (-\psi_{t,k} + \varphi_t + \Lambda)\le C_1^{(n+1)/n} (-\psi_{t,k} + \Lambda) $$ where the second inequality follows from the normalization $\sup_X \varphi_t = 0$. Multiply both sides of above by suitable $0<\alpha = \alpha(\chi, \omega_X)>0$ such that $C_1^{(n+1)/n}  \alpha  $ is less than the alpha invariant of the K\"ahler manifold $(X, (C_2 + 1) \omega_X)$ where $\chi \le C_2 \omega_X$ for some $C_2>0$, and integrate the resulted inequality over $\Omega_s$. We thus obtain
\begin{equation}\label{eqn:trudigner}
\int_{\Omega_s} \exp\Big( \alpha \frac{(v-s)^{(n+1)/n}}{A_{s,k}^{1/n}}   \Big) \omega_X^n \le C \int_{X} \exp\Big( - C_1^{(n+1)/n}  \alpha \psi_{t,k} \Big) \omega_X^n \le C_3,
\end{equation}
for some uniform constant $C_3>0$ independent of $t$. In this last inequality we apply the $\alpha$-invariant estimate for quasi-PSH functions on compact K\"ahler manifolds \cite{H, Ti}. By Young's inequality, \neweqref{eqn:trudigner}  implies that for some $C_4>0$ depending additionally on $p>n$ and ${\mathrm{Ent}}_p(\omega_t)$  (cf. \cite{GPT})
\begin{equation}\label{eqn:key 3}
\int_{\Omega_s} (v-s)^{(n+1)p/n} e^F \omega_X^n \le C_4 A_{s,k}^{p/n} \to C_4 A_{s}^{p/n},
\end{equation}
as $k\to\infty$, noting that the left-hand side of \neweqref{eqn:key 3} is independent of $k$. On the other hand, by H\"older inequality and \neweqref{eqn:key 3} we have
$$
A_s \le \frac{1}{V}\Big(\int_{\Omega_s} (v-s)^{\frac{p(n+1)}{n}} e^F \omega_X^n \Big)^{\frac{n}{p (n+1)}} \Big(\int_{\Omega_s} e^F \omega_X^n \Big)^{\frac{1}{p'}}\le C_5 A_s^{\frac{1}{n+1}} \Big(\int_{\Omega_s} e^F \omega_X^n \Big)^{\frac{1}{p'}},
$$
where $p'>1$ satisfies $\frac{n}{p(n+1)} + \frac{1}{p'} = 1$. This implies that $A_s\le C_6 \Big(\int_{\Omega_s} e^F \omega_X^n \Big)^{\frac{1+n}{n p'}}$, i.e.
\begin{equation}\label{eqn:key 4}
\int_{\Omega_s} (v-s)e^F \omega_X^n \le C_6  \Big(\int_{\Omega_s} e^F \omega_X^n \Big)^{1+\delta_0},
\end{equation}
for $1+\delta_0 = \frac{1+n}{n p'}$ with $\delta_0 = \frac{p-n}{np}>0$. We denote $\phi(s) = \int_{\Omega_s}e^F \omega_X^n$. Then \neweqref{eqn:key 4} yields easily that 
\begin{equation}\label{eqn:key 5} r\phi(s+r) \le C_6 \phi(s)^{1+\delta_0},\quad \forall s\ge 0, \, r>0\end{equation}

\smallskip

\noindent{\bf Step 4.} 
By the assumption $\| v\|_{L^1(X,\omega_t^n)}\le V_0$ we have
$$\int_{\Omega_0} v e^F \omega_X^n \le \frac{1}{c_t} \int_X |v| \omega_t^n\le \frac{V_0 V}{V_t}\le V, $$
and this implies that for any $s>0$
\begin{equation}\label{eqn:key 6}
\phi(s) = \int_{\Omega_s} e^F \omega_X^n \le \frac{1}{s}\int_{\Omega_0}v e^F \omega_X^n \le \frac{V}{s}.
\end{equation}
So we can pick  $s_0 =  (2C_6)^{1/\delta_0} V$ to guarantee that $\phi(s_0)^{\delta_0} < 1/2C_6$. Given \neweqref{eqn:key 5}, we can apply the De Giorgi type iteration argument of Kolodziej \cite{K} to conclude that $\phi(s) = 0$ for any $s> S_\infty$ with 
$$S_\infty = s_0 + \frac{1}{1-2^{-\delta_0}}=(2C_6)^{1/\delta_0} V+\frac{1}{1-2^{-\delta_0}} .$$
This means that $v\le S_\infty$ and we finish the proof of the lemma.

\bigskip


\noindent{\em Proof of Theorem \ref{thm:main1}. }~ Fix a point $x\in X$. We let $v(y) = - G_t(x,y)$ be the Green's function of $\omega_t$ centered at $x$. This $v$ satisfies the assumptions in Lemma \ref{lemma key}, i.e. $v\in L^1(X,\omega_t^n)$, $v$ is smooth on $X\backslash \{x\}$ and 
$$\int_X v( y) \omega_t^n(y) = 0,\quad \Delta_{\omega_t} v(y) = -\frac{1}{V_t} \ge -\frac{1}{V_0}\mbox{ for } y\in \{v\ge 0\}$$
Lemma \ref{lemma key} gives a constant $C>0$ depending only on $n, p, \omega_X, \chi, N$ such that 
$$v\le C(1+ \| v\|_{L^1(X,\omega_t^n)}) $$
which implies that $\inf_{y\in X}G_t(x,y)\ge - C(1+ \|G_t(x,\cdot)\|_{L^1(X,\omega_t^n)} )$.

\smallskip

On the other hand, since $-C_l = \inf_{y\in X}G_t(x,y)$ is a lower bound of $G_t(x,\cdot)$, we have
$$\| G_t(x,\cdot)\|_{L^1(X,\omega_t^n)}\le \int_X |G_t(x,\cdot) + C_l|\omega_t^n + C_l V_t\le 2 C_l V_t.$$

\section{Proof of Theorem \ref{thm:main2}}\setcounter{equation}{0}
Given the parameters $\epsilon>0$, $N>0$ and $\gamma\in (0,1)$, we fix a K\"ahler  metric $\omega_t\in {\mathcal M}_t'(N,\epsilon,\gamma) \cup  {\mathcal M}_t''(N,\epsilon,\gamma) \cup \tilde {\mathcal M}''(N,\epsilon,\gamma)$. We will denote $G_t$ the associated Green's function of $(X,\omega_t)$. As in the last section, we let $\omega_t = \chi +  t\omega_X + \ddbar \varphi_t$ be the solution to the following
complex Monge-Amp\`ere equation
\begin{equation}\label{eqn:MA se3}
(\chi + t \omega_X + \ddbar \varphi_t)^n = c_t e^{F_{\omega_t}} \omega_X^n,\quad \sup_X\varphi_t = 0.
\end{equation}

\smallskip

In this section, we say a constant $C>0$ is {\em uniform} if it depends only on $n, \chi, \omega_X$ and the given parameters $\epsilon, N, \gamma$.

\smallskip

Since $\chi$ is assumed to be nonnegative in Theorem \ref{thm:main2}, by Lemma \ref{lemma GPTW}, there is a uniform constant $C_0>0$ such that for each $t\in (0,1]$ and $\varphi_t$ satisfying \neweqref{eqn:MA se3}
\begin{equation}
\label{eqn:new C0}
\sup_{X} |\varphi_t| \le C_0
\end{equation}

We observe the following estimate on $L^2$-norm of $\nabla \varphi$.

\begin{lemma}\label{lemma gradient new}
Suppose $\omega = \omega_X + \ddbar \varphi$ is a metric such that $e^F = {\omega^n}/{\omega_X^n}\in L^{1+\epsilon}(X,\omega_X^n)$, then 
\begin{equation}\label{eqn:gradient L2}
\int_X |\nabla \varphi|_{\omega_X}^2 \omega_X^n\le C,
\end{equation}
for some $C>0$ depending only on $n,\epsilon, \omega_X$ and $\| e^F\|_{L^{1+\epsilon}}$.
\end{lemma}
\noindent{\em Proof. } If we normalize $\varphi$ such that $\sup_X \varphi = 0$, by the $L^\infty$-estimates of Kolodziej \cite{K, GPT}, we have $\sup_X|\varphi|\le C$ for some $C$ depending on $n, \omega_X$ and $\| e^F\|_{L^{1+\epsilon}}$. Then we calculate
\bea
\int_X(-\varphi) (e^F - 1) \omega_X^n & \nonumber= & \int_X (-\varphi) (\omega^n  - \omega_X^n)\\
& \nonumber= & \int_X i\partial \varphi\wedge \bar \partial \varphi\wedge (\omega^{n-1} + \cdots + \omega_X^{n-1})\\
&\ge \nonumber&\frac{1}{n} \int_X |\nabla \varphi|^2_{\omega_X} \omega_X^n .
\eea
The lemma follows straightforwardly from this.

\smallskip

From Lemma \ref{lemma key}, we easily get that
\begin{lemma}\label{lemma 1}
Suppose $v\in C^2(X)$ satisfies \begin{equation}\label{eqn:s3 1}|\Delta_{\omega_t} v|\le 1,\quad\mbox{ and } \int_X v \omega_t^n = 0,\end{equation} then there is a uniform constant $C>0$ such that
$$\sup_X |v|\le C( 1+  \| v\|_{L^1(X,\omega_t^n)}).$$
\end{lemma}

We further have the following lemma which asserts that the function $v$ satisfying \neweqref{eqn:s3 1} is in fact bounded uniformly in $L^\infty$-norm.
\begin{lemma}\label{lemma 2}
Under the same assumptions as in Lemma \ref{lemma 1}, we have $\| v\|_{L^1(X,\omega_t^n)}\le C$ for some uniform constant $C>0$, in particular by Lemma \ref{lemma 1} this implies $\sup_X |v|\le C$.
\end{lemma}
Assuming Lemma \ref{lemma 2}, we see how it yields the $L^1$-bound on the Green's function $G_t$ of $(X,\omega_t)$. For any fixed point $x\in X$, we view the Green's function $G_t(x,\cdot) = G_t(x,y)$ as a function of $y$. Consider the equation
\begin{equation}\label{eqn:s3 2}\Delta_{\omega_t}  v = - \chi_{\{G_t\le 0\}} + \frac{1}{V_t}\int_{\{G_t\le 0\}} \omega_t^n,\quad\mbox{and } \int_X v \omega_t^n = 0.\end{equation}
Take a sequence of {\em smooth} and {\em bounded} functions $f_k$ that converge pointwise (in fact uniformly) to the bounded function $ - \chi_{\{G_t\le 0\}} + \frac{1}{V_t}\int_{\{G_t\le 0\}} \omega_t^n$ and satisfies $\int_X f_k\omega_t^n = 0$. Let $v_k$ be the smooth solution to $\Delta_{\omega_t} v_k = f_k$ with $\int_X v_k \omega_t^n = 0$. It follows from standard elliptic theory that for fixed $t>0$, $v_k$ converges uniformly to $v$, which is the solution to \neweqref{eqn:s3 2}. By Green's formula, we have
$$v_k(x) = \frac{1}{V_t}\int_X v_k \omega_t^n + \int_X G_t(x,\cdot) (-\Delta_{\omega_t} v_k) \omega_t^n = \int_X G_t(x,\cdot) (-f_k) \omega_t^n.$$
From Lemma \ref{lemma 2} we have $|v_k(x)|\le C$, for each $k$, since $|f_k|\le 2$, say. Letting $k\to \infty$ and  we get by the choice of $f_k$ and the normalization $\int_X G_t(x,\cdot)\omega_t^n = 0$ that 
$$\Big| \int_{\{G_t(x,\cdot)\le 0\}} G_t(x,\cdot)\omega_t^n \Big|\le C.$$
Since $|\int_{\{G_t(x,\cdot)\le 0\}} G_t(x,\cdot)\omega_t^n |= \int_{\{G_t(x,\cdot)\ge 0\}} G_t(x,\cdot)\omega_t^n$, this easily gives the $L^1(X,\omega_t^n)$-bound of $G_t(x,\cdot)$. This finishes the proof of (i) in Theorem \ref{thm:main2}, assuming Lemma \ref{lemma 2}.

\bigskip

Now we turn to the proof of Lemma \ref{lemma 2}. We argue by contradiction. Suppose there is a sequence of K\"ahler metrics $\omega_j = \omega_{t_j}\in {\mathcal M}_t'(N,\epsilon,\gamma) \cup  {\mathcal M}_t''(N,\epsilon,\gamma) \cup \tilde {\mathcal M}''(N,\epsilon,\gamma)$
and
$F_j = F_{\omega_j}$ as defined in \neweqref{eqn:relative volume},  and a sequence of $C^2(X)$ functions $v_j$ satisfying
$$\Delta_{\omega_j} v_j = h_j,\quad \mbox{and }\int_X v_j \omega_j^n = 0,$$
for some function $h_j$ with $\sup_X|h_j|\le 1$, for which Lemma \ref{lemma 2} fails, i.e. they satisfy 
\begin{equation}\label{eqn:0}\| v_j\|_{L^1(X,\omega_j^n)} = \int_X |v_j| \omega_j^n \to \infty \mbox{ ~~as }j\to\infty.\end{equation}
We normalize each $v_j$ by $$\hat v_j = \frac{v_j}{\| v_j\|_{L^1(X, \omega_j^n)}},\quad\mbox{so } \| \hat v_j\|_{L^1(X,\omega_j^n)} = 1.$$
It is clear that $\hat v_j\in C^2(X)$ satisfies 
\begin{equation}\label{eqn:1}\Delta_{\omega_j} \hat v_j  = \frac{h_j}{\| v_j\|_{L^1(X, \omega_j^n)}},\quad \mbox{and } \int_X \hat v_j \omega_j^n = 0.\end{equation}
We can apply Lemma \ref{lemma 1} to conclude that
\begin{equation}\label{eqn:2}
\sup_X |\hat v_j|\le C(1+ \| \hat v_j\|_{L^1(X,\omega_j^n)})\le C,
\end{equation}
for some uniform constant $C>0$ independent of $j$. Multiplying both sides of \neweqref{eqn:1} by $\hat v_j$ and applying integration by parts, we get
\begin{equation}\label{eqn:3}
\int_X |\nabla \hat v_j|^2_{\omega_j} \omega_j^n = \int_X - \frac{h_j \hat v_j}{\| v_j\|_{L^1(X,\omega_j^n)}}\omega_j^n\le \frac{1}{\| v_j\|_{L^1(X,\omega_j^n)}}\to 0,
\end{equation}as $j\to\infty$ by the hypothesis \neweqref{eqn:0}. By the H\"older inequality we have
\bea
\int_X |\nabla \hat v_j|_{\omega_X}\omega_X^n &\nonumber \le &\int_X (|\nabla \hat v_j|_{\omega_j}^2 \tr_{\omega_X} \omega_j)^{1/2}\omega_X^n\\
&\le & \label{eqn:4} \bk{\int_X |\nabla \hat v_j|_{\omega_j}^2 e^{F_j} \omega_X^n}^{1/2} \bk{\int_X( \tr_{\omega_X}\omega_j) e^{-F_j}\omega_X^n}^{1/2}.
\eea
Note that the first factor in \neweqref{eqn:4} satisfies
$$\int_X |\nabla \hat v_j|_{\omega_j}^2 e^{F_j}\omega_X^n = c_{t_j}^{-1} \int_X |\nabla \hat v_j|_{\omega_j}^2 \omega_j^n\to 0$$
by \neweqref{eqn:3}. For the second factor in \neweqref{eqn:4}, we consider different cases of $\omega_j$. 

\smallskip

\noindent (a) If $\omega_j\in {\mathcal M}'_{t_j}(N,\epsilon,\gamma)$, then it holds that
$$\int_X( \tr_{\omega_X}\omega_j) e^{-F_j}\omega_X^n\le \frac{1}{n\gamma} \int_X\omega_j\wedge \omega_X^{n-1}\le C.$$
\smallskip

\noindent (b) If $\omega_j\in {\mathcal M}''_{t_j}(N,\epsilon,\gamma)$, then we have
\bea
\nonumber 
\int_X( \tr_{\omega_X}\omega_j) e^{-F_j}\omega_X^n& = & \frac 1 n\int_X e^{-F_j} (\chi + t_j \omega_X + \ddbar \varphi_{t_j})\wedge \omega_X^{n-1}\\
&\nonumber \le &C \int_X e^{-F_j} \omega_X^n + \int_X|\varphi_{t_j}| |\Delta_{\omega_X} e^{-F_j}| \omega_X^n\\
&\le&\nonumber C,
\eea
by the definition of $ {\mathcal M}''_{t_j}(N,\epsilon,\gamma)$ and \neweqref{eqn:new C0}.

\smallskip

\noindent (c) If $\omega_j\in \tilde {\mathcal M}''(N,\epsilon,\gamma)$, then by a similar calculation we have
\bea
\nonumber 
\int_X( \tr_{\omega_X}\omega_j) e^{-F_j}\omega_X^n& = & \frac 1 n\int_X e^{-F_j} ( \omega_X + \ddbar \varphi_{j})\wedge \omega_X^{n-1}\\
&\nonumber \le &C \int_X e^{-F_j} \omega_X^n + \int_X|\nabla \varphi_{j}|_{\omega_X} |\nabla e^{-F_j}|_{\omega_X} \omega_X^n\\
&\nonumber \le &C \int_X e^{-F_j} \omega_X^n +\bk{ \int_X|\nabla \varphi_{j}|_{\omega_X}^2 \omega_X^n}^{1/2}\bk{\int_X |\nabla e^{-F_j}|^2_{\omega_X} \omega_X^n}^{1/2}\\
&\le&\nonumber C,
\eea
by Lemma \ref{lemma gradient new} and the definition of the set $ \tilde {\mathcal M}''(N,\epsilon,\gamma)$. 

\smallskip

\noindent (d) If $n = 1$, we observe that from \neweqref{eqn:3}, as $j\to \infty$
\begin{equation}\label{eqn:n=1}\int_X |\nabla \hat v_j|^2_{\omega_X} \omega_X= \int_X i\partial \hat v_j \wedge \bar \partial \hat v_j = \int_X |\nabla \hat v_j|^2_{\omega_j} \omega_j \to 0. \end{equation}

\smallskip

Combining all cases discussed above,
%
%
 \neweqref{eqn:4} or \neweqref{eqn:n=1} entail that
\begin{equation}\label{eqn:5}
\int_X |\nabla \hat v_j|_{\omega_X}\omega_X^n \to 0.
\end{equation}
From \neweqref{eqn:2} and \neweqref{eqn:5} we see the the sequence of functions $\{\hat v_j\}$ is uniformly bounded in the Sobolev space $W^{1,1}(X,\omega_X)$ (under the {\bf fixed} metric $\omega_X$). By the Sobolev embedding theorem, there is an embedding
$$W^{1,1}(X,\omega_X) \hookrightarrow L^{q}(X,\omega_X^n)$$
which is compact for any $1\le q< \frac{2n}{2n - 1}$. Therefore, up to a subsequence we have $\hat v_j\to \hat v_\infty$ in $L^{q}(X,\omega_X^n)$. In particular $\hat v_j$ also converge to $\hat v_\infty$ in $L^1(X,\omega_X^n)$ and  in the a.e. sense up to a further subsequence if necessary. 

{We now claim that $\hat v_\infty$  is constant in the a.e. sense. Indeed, for any fixed $C^2$ function $\rho$ on $X$, we have
$$
\Big|\int_X \hat v_j \Delta_{\omega_X} \rho \omega_X^n\Big| = \Big|\int_X \langle\nabla \hat v_j, \nabla  \rho\rangle_{\omega_X} \omega_X^n\Big|\le \|\nabla \rho\|_{\infty} \int_X |\nabla \hat v_j|_{\omega_X} \omega_X^n\to 0
$$
as $j\to\infty$. By the dominated convergence theorem and $\hat v_j\to \hat v_\infty$ a.e. we conclude that
$$\int_X \hat v_\infty \Delta_{\omega_X} \rho \omega_X^n = 0$$ which holds for {\em any} $\rho\in C^2(X)$. By Weyl's lemma, this implies that $\hat v_\infty$ is $\Delta_{\omega_X}$-harmonic, 
}
hence $\hat v_\infty = \alpha_0$  in the a.e. sense for some constant $\alpha_0\in\mathbb R$. We next claim that $\alpha_0\neq 0$. Indeed, from the normalization $\| \hat v_j\|_{L^1(X, \omega_j^n)} = 1$ and \neweqref{eqn:2}, we get
\bea
1 & = &\nonumber c_{t_j}\int_X |\hat v_j| e^{F_j} \omega_X^n\le C \int_X |\hat v_j|^\eta e^{F_j} \omega_X^n\\
 &\le&  C\bk{\int_X |\hat v_j|^{\eta \frac{1+\epsilon}{\epsilon}} \omega_X^n   }^{\epsilon/(1+\epsilon)} \bk{\int_X e^{(1+\epsilon) F_j}\omega_X^n}^{1/(1+\epsilon)} \label{eqn:8}
\eea
and here we take $\eta = \epsilon/(1+\epsilon)$. From \neweqref{eqn:8} we obtain
\begin{equation}\label{eqn:9}
\int_X |\hat v_j| \omega_X^n \ge c_0>0
\end{equation}
for some uniform constant $c_0>0 $. Taking limit and applying the dominated convergence theorem again, we get $\int_X |\hat v_\infty|\omega_X^n \ge c_0>0$, and this implies that $|\alpha_0|>0$.

\medskip

However, this will contradict the second equation in \neweqref{eqn:1}. To see this, we assert that $\lim_{j\to \infty} \int_X \hat v_\infty \omega_j^n = 0$. In fact, by \neweqref{eqn:1} we have
\bea
| \int_X \hat v_\infty \omega_j^n| &= \nonumber& |  \int_X (\hat v_j - \hat v_\infty) \omega_j^n  |\\
&\le &\nonumber c_{t_j}  \int_X | \hat v_j - \hat v_\infty|e^{F_j} \omega_X^n\\
&\le & \nonumber C \int_X | \hat v_j - \hat v_\infty|^\eta e^{F_j} \omega_X^n \quad \mbox{ here }\eta = \frac{\epsilon}{1+\epsilon}\\
&\le &\nonumber C \bk{ \int_X | \hat v_j - \hat v_\infty|\omega_X^n}^{\epsilon/(1+\epsilon)}\bk{\int_X  e^{(1+\epsilon)F_j} \omega_X^n}^{1/(1+\epsilon)}\\ 
&\le & \nonumber C \bk{ \int_X | \hat v_j - \hat v_\infty|\omega_X^n}^{\epsilon/(1+\epsilon)}\to 0
\eea
since $\hat v_j \to \hat v_\infty$ in $L^1(X,\omega_X^n)$. But this is absurd since $ \int_X \hat v_\infty \omega_j^n = \alpha_0 V_{t_j}$ which is strictly away from zero.
This finishes the proof of Lemma \ref{lemma 2}.

\bigskip


Once the $L^1(X,\omega_t^n)$-norm of $G_t$ is achieved, Theorem \ref{thm:main1} provides a lower bound of $G_t$. \neweqref{eqn:Green lower bound} is thus proved. Let $C_l>0$ be the constant in \neweqref{eqn:Green lower bound}, %
i.e. $G_t\ge -C_l$. For notational simplicity, we will denote the {\em positive Green's function} 
\begin{equation}\label{eqn:positive Green}{\mathcal G}_t (x,\cdot)= G_t(x,\cdot) + C_l + 1\ge 1.\end{equation}



We are ready to prove (ii) in Theorem \ref{thm:main2}. To begin with, we show that the $L^q(X,\omega_t^n)$-norm of $G_t$ is uniformly bounded, for any $q<\frac{2n}{2n-2}$, which is optimal in view of the asymptotic behavior of $G_t$ in \neweqref{eqn:G asymp}.

\begin{lemma}\label{lemma 2.3}
For any $q\in (1, \frac{2n}{2n -2})$, there is a uniform constant $C>0$ depending on $q$ such that $G_t$ satisfies
\begin{equation}\label{eqn:green power}
\int_X |G_t(x,y)|^q \omega_t^n(y) \le C, \quad \forall~ x\in X. 
\end{equation}
\end{lemma}

\smallskip
\noindent{\em Proof. }
We break the proof into two steps. We fix a point $x\in X$ and consider the Green's function $G_t(x,y)$ as a function of $y$. It suffices to show \neweqref{eqn:green power} for ${\mathcal G}_t(x,\cdot)$, since $C_l>0$ in \neweqref{eqn:positive Green} is uniform. The first step is to show the $L^q(X,\omega_t^n)$ bound of ${\mathcal G}_t(x,\cdot)$ for any $q< 1+1/n$, then we can apply an iteration argument similar to the Moser iteration process to improve the exponent $q$.

\smallskip

\noindent{\bf Step 1.} We will show the $L^{1+ \frac{1}{r_0}}(X, \omega_t^n)$-norm of ${\mathcal G}_t(x,\cdot)$ is uniformly bounded, for {\em any} $r_0>n$. The argument is based on the $L^1$-bound of ${\mathcal G}_t(x,\cdot)$ in (i) of Theorem \ref{thm:main2}.

\smallskip

Fix a large $k\gg1$ and consider the function $H_k(y) = \min\{{\mathcal G}_t(x,y), k\}$. By smoothing $H_k$ if necessary we may assume it is a smooth function and $H_k(y)$ converges monotonically increasingly to ${\mathcal G}_t(x,y)$ as $k\to \infty$. We solve the following equation{\small
\begin{equation}\label{eqn:new key 1}
\left\{\begin{array}{ll}
&\Delta_{\omega_t} u_k = - H_k^{1/r_0} + \frac{1}{V_t} \int_X H_k^{1/r_0} \omega_t^n,\\
&\frac{1}{V_t} \int_X u_k\omega_t^n = 0.\end{array}\right.
\end{equation}
}%
Equation \neweqref{eqn:new key 1} admits a unique smooth solution since the smooth function on the right-hand side of the first equation has integral $0$. To deal with the unbounded term $- H_k^{1/r_0}$ on the right-hand side of \neweqref{eqn:new key 1} and non-uniform ellipticity of the linear operator $\Delta_{\omega_t}$,  we  consider again an auxiliary complex Monge-Amp\`ere equation

\begin{equation}\label{eqn:aMA new}
(\chi + t\omega_X + \ddbar \psi_k)^n = \frac{H_k^{n/r_0} + 1}{V_t^{-1} \int_X (H_k^{n/r_0} + 1) \omega_t^n} \omega_t^n = c_t \frac{H_k^{n/r_0} + 1}{B_k} e^{F_{\omega_t}} \omega_X^n,
\end{equation}
with $\sup_X \psi_k = 0$ and $B_k = \int_X (H_k^{n/r_0} + 1) e^{F_{\omega_t}} \omega_X^n$. We stress that this auxiliary Monge-Amp\`ere equation plays a very different role from the auxiliary Monge-Amp\`ere equation introduced in the proof of Lemma 2.
We note that 
\begin{equation}\label{eqn:Bk}
V\le B_k\le V + \Big( \int_X H_k e^{F_{\omega_t}} \omega_X^n \Big)^{n/r_0} \Big(\int_X e^{F_{\omega_t}} \omega_X^n \Big)^{(r_0 - n)/r_0} \le C(V),
\end{equation}
and the upper bound holds because of $0<H_k\le {\mathcal G}_t$ for any $k$ and the integral bound of ${\mathcal G}_t$ in (i) of Theorem \ref{thm:main2}. We note that the $p$-th (for some $p>n$) entropy of  the function on the right-hand side of \neweqref{eqn:aMA new} satisfies
\bea
\nonumber && \frac{1}{B_k}\int_X (H_k^{n/r_0} + 1  )  \Big| -\log B_k + F_{\omega_t} + \log ( 1+ H_k^{n/r_0})   \Big|^p e^{F_{\omega_t}} \omega_X^n\\
\label{eqn:newest 1} &\le & \frac{|\log B_k|^p}{B_k} \int_X (H_k^{n/r_0} + 1  ) e^{F_{\omega_t}} \omega_X^n + \frac{1}{B_k} \int_X (H_k^{n/r_0} + 1  ) [\log(H_k^{n/r_0} + 1  )]^p e^{F_{\omega_t}} \omega_X^n \\
\nonumber & & \quad  + \frac{1}{B_k} \int_X (H_k^{n/r_0} + 1  ) |F_{\omega_t}|^p e^{F_{\omega_t}} \omega_X^n
\le C
\eea
for some uniform constant $C>0$ depending on $n, p,\epsilon, \chi,\omega_X$ and $\| e^{F_{\omega_t}}\|_{L^{1+\epsilon}(\omega_X^n)}$. Here the first term in \neweqref{eqn:newest 1} is bounded due to the estimate of the constant $B_k$ in \neweqref{eqn:Bk} and H\"older inequality along with the uniform $L^1(X,\omega_t^n)$-bound of $H_k$; the second term is bounded because of $\log(1+ x)\le C_\delta x^\delta$ for any $\delta>0$ and the $L^1(X,\omega_t^n)$-bound of $H_k$; and the last term is bounded again by H\"older inequality. We can now apply Lemma \ref{lemma GPTW} to conclude that
$$\sup_X | \psi_k |\le C,$$
We know $\varphi_t$ also satisfies \neweqref{eqn:new C0}, so it follows that
\begin{equation}\label{eqn:new key 2}
\sup_X |\psi_k - \varphi_t|\le C,
\end{equation} for some uniform constant $C>0$, which is independent of $k$.
We now consider the function
\begin{equation}\label{eqn:def v}v: = (\psi_k - \varphi_t) - \frac{1}{V_t} \int_X (\psi_k - \varphi_t) \omega_t^n + \varepsilon' u_k ,
\end{equation}
where $\varepsilon'>0$ is a suitable constant to be chosen later. By definition it follows $\frac{1}{V_t}\int_X v \omega_t^n = 0$ and $v$ is a smooth function. We then calculate the Laplacian of $v$ in \neweqref{eqn:def v}
\bea
\Delta_{\omega_t} v \nonumber &= & \tr_{\omega_t} \omega_{t, \psi_k} - n + \varepsilon '  \Delta_{\omega_t} u_k\\
\nonumber &\ge & n \Big( \frac{\omega_{t,\psi_k}^n}{\omega_t^n} \Big)^{1/n} - n - \varepsilon' H_k^{1/r_0} + \frac{\varepsilon'}{V_t} \int_X H_k^{1/r_0} \omega_t^n\\
\nonumber &= & n B_k^{-1/n} ( H_k^{n/r_0} + 1 )^{1/n} - n - \varepsilon' H_k^{1/r_0} + \frac{\varepsilon'}{V_t} \int_X H_k^{1/r_0} \omega_t^n\\
\nonumber &\ge & n B_k^{-1/n}  H_k^{1/r_0} - n - \varepsilon' H_k^{1/r_0} \ge -n,
\eea
if we choose $\varepsilon' = n C(V)^{-1/n}$, where $C(V)$ is the upper bound of $B_k$ in \neweqref{eqn:Bk}. We apply the Green's formula to the function $v$ at $x$
\bea
v(x) & = & \frac{1}{V_t} \int_X v \omega_t^n + \int_X G_t(x,\cdot) (-\Delta_{\omega_t} v) \omega_t^n =  \int_X {\mathcal G}_t(x,\cdot) (-\Delta_{\omega_t} v) \omega_t^n \nonumber\\
&\le & \int_{X} {\mathcal G}_t(x,\cdot) n \omega_t^n\le C,\nonumber
\eea
where the last inequality follows from the uniform $L^1(X,\omega_t^n)$-bound of ${\mathcal G}_t(x,\cdot)$, as proved in (i) of Theorem \ref{thm:main2}. It then follows from \neweqref{eqn:new key 2} and \neweqref{eqn:def v} that $u_k(x)\le C$ for a uniform constant $C>0$.


\smallskip


We apply again the Green's formula to the function $u_k$ at $x$ to get
\bea
C\ge u_k(x)\nonumber & = &\frac{1}{V_t} \int_X u_k \omega_t^n + \int_X {\mathcal G}_t(x,\cdot) (-\Delta_{\omega_t} u_k)\omega_t^n\\
\nonumber &= & \int_X {\mathcal G}_t(x,y) \Big( H_k(y)^{1/r_0} - \frac{1}{V_t} \int_X H_k^{1/r_0} \omega_t^n \Big)\omega_t^n(y),
\eea
which yields that
\begin{equation}\label{eqn:newest 2}\int_X {\mathcal G}_t(x,y) H_k(y)^{1/r_0} \omega_t^n(y)\le C + C \int_X H_k^{1/r_0} \omega_t^n \le C
\end{equation}
where we apply the H\"older inequality and the $L^1(X,\omega_t^n)$-bound of ${\mathcal G}_t(x,\cdot)$ and $H_k$. Finally letting $k\to\infty$ we conclude from monotone convergence theorem and \neweqref{eqn:newest 2} that
\begin{equation}\label{eqn:step 1}
 \int_X {\mathcal G}_t(x,y)^{1+ \frac{1}{r_0}} \omega_t^n(y)\le C(r_0).
\end{equation}
Here we write the uniform constant as $C(r_0)$ to emphasize its addintional dependence on $r_0>n$.

\medskip

\noindent{\bf Step 2.} We now iterate the arguments in {\bf Step 1}, given the $L^{1+1/r_0}(X,\omega_t^n)$-bound \neweqref{eqn:step 1} of ${\mathcal G}_t$ for any $r_0>n$. We take an arbitrary $r_1\in (0,r_0)$ such that $\frac{n}{r_1} < 1+ \frac{1}{r_0}$, and replace $r_0$ by $r_1$ in the equations \neweqref{eqn:new key 1} and \neweqref{eqn:aMA new}. Using the better $L^{1+1/r_0}(X,\omega_t^n)$-bound \neweqref{eqn:step 1} of ${\mathcal G}_t$ instead of just the $L^1(X,\omega_t^n)$-bound as in {\bf Step 1},  we can repeat the arguments above to obtain the $L^{1+1/r_1}(X,\omega_t^n)$-bound of ${\mathcal G}_t$. Iterating this process we eventually get for {\em any} $r_l\in (0, r_{l-1})$ with $\frac{n}{r_l}< 1+ \frac{1}{r_{l-1}}$, the $L^{1+1/r_l}(X,\omega_t^n)$-bound of ${\mathcal G}_t$ is achieved, which depends in addition on $l$. A simple combinatorial argument shows that one can take any $r_l>0$ such that 
$$\frac{1}{r_l} < \frac{1}{n-1} - \frac{1}{n^l}\frac{1}{n(n-1)}.$$
For $l$ large enough, this implies the bound on $L^q(X,\omega_t^n)$-norm of ${\mathcal G}_t(x,\cdot)$ for any $q< 1+ \frac{1}{n-1} = \frac{n}{n-1}$. As we mentioned at the beginning, the $L^q(X,\omega_t^n)$-bound of $ G_t(x,\cdot)$ follows from the definition of ${\mathcal G}_t(x,\cdot)$ in \neweqref{eqn:positive Green}.
\bigskip

Now we derive the $L^q(X,\omega_t^n)$ bound on the gradient of ${\mathcal G}_t(x,\cdot)$ for any $q\in (1, \frac{2n}{2n-1})$. First we observe the following elementary estimate which follows easily from the Green's formula.
\begin{lemma}\label{lemma gradient}
For any $\beta>0$ we have
\begin{equation}\label{eqn:grad 1}
\int_X \frac{ |\nabla_y {\mathcal G}_t(x,y)|_{\omega_t(y)}^2   }{{\mathcal G}_t(x,y)^{1+\beta}} \omega_t^n(y) \le 1/\beta.
\end{equation} 
\end{lemma}
\noindent{\em Proof. }
The function $u(y): = {\mathcal G}_t(x,y)^{-\beta}$ is a continuous function with $u(x) = 0$ and $u\in C^\infty(X\backslash \{x\})$. By the Green's formula we have
\bea0= u(x) \nonumber &= & \frac{1}{V_t} \int_X u \omega_t^n + \int_X {\mathcal G}_t(x, \cdot) (-\Delta_{\omega_t} u) \omega_t^n\\
\nonumber &= & \frac{1}{V_t} \int_X u \omega_t^n -\beta \int_X \frac{|\nabla {\mathcal G}_t(x,\cdot)|^2_{\omega_t}}{ {\mathcal G}_t(x,\cdot)^{1+\beta}} \omega_t^n.
\eea
Here we have applied the integration by parts, which can be justified by the asymptotic behavior of ${\mathcal G}_t(x,y)$ as in \neweqref{eqn:G asymp}. The lemma then follows easily from the fact that $0\le u \le 1$.

\begin{lemma}
\label{lemma mix gradient}
For any $\delta\in (0, 2/n)$ and $\beta>0$, there is a uniform constant $C>0$ depending on $\delta$ and $\beta$ such that 
\begin{equation}\label{eqn:grad 2}
\int_X {\mathcal G}_t(x, y) ^{1- \frac{(1+\beta) (2-n\delta)}{2n} } |\nabla G_t(x, y)|^{\frac 2n -\delta}_{\omega_t(y)} \omega_t(y)\le C.
\end{equation}

\end{lemma}
\noindent{\em Proof. }
Given Lemma \ref{lemma gradient}, the proof of Lemma \ref{lemma mix gradient} is essentially the same as that in {\bf Step 1} of Lemma \ref{lemma 2.3}. We use similar notations as in Lemma \ref{lemma 2.3}. Let
$$H_k(y) =\widetilde{ \min}\Big\{  \frac{ |\nabla_y {\mathcal G}_t(x,y)|_{\omega_t(y)}^{2}   }{{\mathcal G}_t(x,y)^{(1+\beta)}}  , k\Big\}$$
where $\widetilde{\min}$ denotes a suitable smoothing of the $\min$ function. We can solve equations \neweqref{eqn:new key 1} and \neweqref{eqn:aMA new} with this $H_k$ and any $r_0>n$. With the estimate \neweqref{eqn:grad 1} in Lemma \ref{lemma gradient}, the same argument gives
$$\int_X {\mathcal G}_t(x,y) H_k(y)^{1/r_0} \omega_t^n(y) \le C.$$
Letting $k\to\infty$, this  yields 
$$\int_{X} {\mathcal G}_t(x,y)^{1- (1+\beta)/r_0}  |\nabla_y {\mathcal G}_t(x,y)|_{\omega_t(y)}^{2/r_0} \omega_t^n(y)\le C. $$
The lemma follows by setting $\delta = \frac{2}{n} - \frac{2}{r_0}$.

\medskip

With Lemmas \ref{lemma gradient} and \ref{lemma mix gradient}, we are ready to derive the $L^s(X,\omega_t^n)$-bound of $\nabla G_t(x,\cdot)$ for  $s\in [1,\frac{2n}{2n-1})$.

\begin{lemma}\label{lemma gradient Lp}
For any given $s\in [1, \frac{2n}{2n-1})$, there is a uniform constant $C>0$ depending on $s$ such that
\begin{equation}\label{eqn:grad lp}
\int_X |\nabla G_t(x, y)|^s_{\omega_t(y)} \omega_t^n(y) \le C.
\end{equation}
\end{lemma}
\noindent{\em Proof. }
We fix a constant $1\le s< \frac{2n}{2n-1}$ and a point $x\in X$. The Green's function $G_t(x,y)$ is viewed as a function of $y$ and all the integrals below are integrated over $y\in X$. Then we calculate
\bea
& \nonumber \int_X  |\nabla G_t(x, y)|^s_{\omega_t(y)} \omega_t^n(y)  = \int_X ( |\nabla G_t|_{\omega_t} ^{s \delta} {\mathcal G}_t^a )  \frac{|\nabla G_t|_{\omega_t}^{(1-\delta) s}  }{{\mathcal G}_t^b} {\mathcal G}_t^{b-a} \omega_t^n\\
\le \label{eqn:three} & \Big(\int_X |\nabla G_t|_{\omega_t}^{s\delta p} {\mathcal G}_t^{ap}\omega_t^n \Big)^{1/p}      \Big( \int_{X} \frac{|\nabla G_t|_{\omega_t}^{(1-\delta) s q}   }{{\mathcal G}_t^{b q}} \omega_t^n \Big)^{1/q}          \Big(  \int_X {\mathcal G}_t^{(b-a)r} \omega_t^n   \Big)^{1/r}
\eea
where the inequality follows from the generalized H\"older inequality and the numbers involved are chosen as follows.
$$\frac{1}{p} + \frac{1}{q} + \frac{1}{r} = 1, \, 0< \delta \ll 1, \, 0<a<b, \, (b-a) r< \frac{n}{n-1}.$$
Furthermore for a small $\beta>0$ we can pick these numbers as 
\begin{equation}\label{eqn:choice}(1 - \delta) s q = 2,\, \, b q = 1+ \beta,\,\, s \delta p = \frac{2}{n} - \beta,\,\, a p = 1 - \frac{(1+\beta) (2 - n \beta)}{2n}. \end{equation}
With these choices of parameters, Lemma \ref{lemma gradient} and Lemma \ref{lemma mix gradient} imply the first two factors in \neweqref{eqn:three} are bounded. It only remains to verify the last factor in \neweqref{eqn:three} is also bounded, and by Lemma \ref{lemma 2.3}, it suffices to make sure that $(b-a) r < \frac{2n}{2n - 2} = \frac{n}{n-1}$. From \neweqref{eqn:choice}, we derive
$$q = \frac{2}{(1-\delta) s},\, p = \frac{{2} - n \beta}{ ns \delta},  \, b = \frac{(1+\beta) (1-\delta) s}{2},\, a = \frac{ n s \delta}{{2} - n \beta} \Big( 1 - \frac{(1+\beta) (2-n\beta)}{2n}   \Big) .$$ Therefore, we have
{
\bea
\nonumber && (b-a) r = : \Psi(\delta)\\
\nonumber &= & \Big( \frac{(1+\beta) (1-\delta) s}{2} - \frac{ n s \delta}{{2} - n \beta} \big( 1 - \frac{(1+\beta) (2-n\beta)}{2n}   \big)  \Big) \frac{1}{1 - \frac{(1-\delta) s}{2} - \frac{n s\delta}{2 - n \beta}  } .
\eea }%
We note that as $s< \frac{2n}{2n-1}$
$$\Psi(0)= \frac{(1+\beta) s}{2} \frac{1}{1- \frac{s}{2}} = (1+\beta) \frac{s}{2 -s}< \frac{2n}{2n-2} $$
if $\beta>0$ is chosen small enough. Since $\Psi(\delta)$ is continuous in $\delta\ge 0$, $\Psi(\delta)< \frac{2n}{2n-2}$ if $\delta>0$ is sufficiently small. This verifies the desired inequality $(b-a)r < \frac{n}{n-1}$. The proof of Lemma \ref{lemma gradient Lp} is completed. 

\bigskip

As a corollary of Lemma \ref{lemma gradient Lp}, we have the following Sobolev-Morrey type inequality for the metric $\omega_t$. 
\begin{corollary}
For any $p>2n$ there is a uniform constant $C>0$ depending on $p$ such that 
$$\sup_X \Big| u - \frac{1}{V_t} \int_X u \omega_t^n  \Big| \le C \Big(\int_X |\nabla u|_{\omega_t}^p \omega_t^n \Big)^{1/p}, \,\, \forall u\in C^1(X).$$
\end{corollary}
\noindent{\em Proof. }
This follows immediately from the Green's formula below, Lemma \ref{lemma gradient Lp} and H\"older inequality
$$u(x) - \frac{1}{V_t} \int_X u \omega_t^n = \int_X \langle \nabla_y G_t(x, y), \nabla u(y) \rangle_{\omega_t(y)} \omega_t^n(y), $$ noting that the conjugate exponent $p^* = \frac{p}{p-1}< \frac{2n}{2n-1}$.%

\medskip

\noindent{\em Proof of (ii) in Theorem \ref{thm:main2}. } This follows from Lemma \ref{lemma 2.3} and Lemma \ref{lemma gradient Lp}.

\bigskip

We conclude this section by comparing the lower bound of the Green's function obtained in Theorem \ref{thm:main2} with the classical one in Cheng-Li \cite{CL}. Let $\omega = \omega_X+ \ddbar \varphi \in [\omega_X]$ be a K\"ahler metric with $\| e^{F_{\omega}}\|_{L^{1+\epsilon}(X,\omega_X^n)}\le N$. Suppose $\ric(\omega)\ge -\kappa' \omega$ for some $\kappa'\ge 0$, then from \cite{GPTW1, FGS} we know ${\mathrm{diam}}(X,\omega)\le C(n, \omega_X, N)$. Then Cheng-Li's estimate \neweqref{eqn:Cheng Li} implies the Green's function associated with $\omega$ is bounded below. 

\smallskip

We show now that, in the K\"ahler setting, under a Kolodziej type \cite{K} condition on the volume form, Theorem \ref{thm:main2} implies the lower bound of the Green's function under the less restrictive assumption of a lower bound of the {\em scalar curvature}:

\begin{corollary}\label{cor CL}
For any K\"ahler metric $\omega\in [\omega_X]$, if its relative volume form $\| e^{F_\omega}\|_{L^{1+\epsilon}(X,\omega_X^n)}\le N$ for some $\epsilon>0$, $N>0$, and its scalar curvature $R(\omega)\ge - \kappa$ for some $\kappa\ge 0$, the Green's function $G$ of $(X,\omega)$ satisfies
$$\inf_{y\in X} G(x, y)\ge -C,\quad \forall x\in X,$$
for some constant $C>0$ depending on $n, \omega_X, \epsilon, N$ and $\kappa$.
\end{corollary}

\noindent{\em Proof. } We claim that under the assumption of a scalar curvature lower bound, the relative volume form satisfies $\inf_X e^{F_\omega}\ge \delta'$ for some $\delta'>0$. Hence it follows that $\omega\in {\mathcal M}'(N,\epsilon,\gamma)$ with $\gamma = 1/\delta'$ and (i) in Theorem \ref{thm:main2}  implies a lower bound on the Green's function $G$. 

\smallskip

To see the claim, note that by definition of the Ricci curvature, $\ric(\omega) = \ric(\omega_X) + \ddbar (- F_\omega)$, so the scalar curvature of $\omega$ satisfies $R(\omega) = \tr_\omega(\ric(\omega_X)) + \Delta_{\omega} (-F_\omega)$. The assumption that $R(\omega)\ge -\kappa$ implies $\Delta_\omega (-F_{\omega})\ge -\kappa - \tr_\omega(\ric(\omega_X))$. By Kolodziej's $L^\infty$ estimate \cite{K} (see also \cite{GPT}), the K\"ahler potential $\varphi$ of $\omega = \omega_X+ \ddbar \varphi$ is bounded, i.e. $\| \varphi\|_{L^\infty}\le C$ for some $C=C(\epsilon, N)>0$, if we normalize $\sup_X\varphi = 0$.
For a constant $A>0$ to be determined, we calculate
\bea \Delta_\omega(- F_\omega - A\varphi)& \ge&\nonumber -\kappa - \tr_\omega(\ric(\omega_X)) + A \tr_\omega \omega_X- An\\
&\ge & \tr_\omega \omega_X -\nonumber \kappa - An\\
&\ge & n \Big(\frac{\omega_X^n}{\omega^n} \Big)^{1/n} - \kappa - An \nonumber = n e^{-F_\omega/n} - \kappa - An
\eea
where we take $A = C' + 1$ and $C'>0$ is an upper bound of the Ricci curvature $\ric(\omega_X)$. Applying maximum principle, we get at the maximum point of $-F_\omega - A\varphi$, $e^{-F_\omega}\le (A+ \kappa/n)^n$. Combined with the $L^\infty$ bound of $\varphi$, this easily shows the upper bound of $-F_\omega$.

\medskip
We observe that typically, the assumption of lower bound for the scalar curvature is much more difficult to work with than the assumption of lower bound for the Ricci curvature. There seem to be far fewer results under this assumption. One of which, in a very different direction, is the recent result of  Munteanu-Wang \cite{MW} on the decay of the Green's function on a real three-dimensional complete manifold with scalar curvature bounded from below.

\section{Applications}\setcounter{equation}{0}
\label{section 4}

In this section we discuss some applications of the estimates on the Green's functions in Theorem \ref{thm:main2}, and provide the proof of Theorems \ref{thm:gradient} and \ref{thm:C2}. Given the parameters $\epsilon, N,$ and $\gamma$, for each $t\in (0,1]$, we fix a K\"ahler metric 
$\omega_t \in  {\mathcal M}_t'(N,\epsilon,\gamma) \cup  {\mathcal M}_t''(N,\epsilon,\gamma) \cup \tilde {\mathcal M}''(N,\epsilon,\gamma)$.

\smallskip

Recall that $\chi$ a closed and {\em nonnegative} $(1,1)$-form such that its class $[\chi]$ is big. By the definition of $F_{\omega_t}$ in \neweqref{eqn:relative volume}, the metric $\omega_t = \chi + t\omega_X + \ddbar \varphi_t$ satisfies the following complex Monge-Amp\`ere equation
\begin{equation}\label{eqn:MA s3}
(\chi + t \omega_X + \ddbar \varphi_t)^n = c_t e^{F_{\omega_t}} \omega_X^n,\quad \mbox{and }\sup_X  \varphi_t = 0,
\end{equation}
where as before $c_t = V_t/V$ is a normalizing constant. By the assumptions on $\omega_t$, the Green's function $G_t$ associated with $\omega_t$ satisfies the estimates stated in Theorem \ref{thm:main2}.


 Since $[\chi]$ is assumed to be big, Kodaira's lemma implies that there is an effective divisor $D$ 
 such that
\begin{equation}\label{eqn:Kodaira}\chi - \varepsilon_0 \ric(h_D) \ge \delta_0 \omega_X,\end{equation}
where $h_D$ is a Hermtian metric on the line bundle $[D]$ associated to $D$ and $\varepsilon_0>0$ and $\delta_0>0$ are fixed constants depending only on $\chi, \omega_X$. Let $s_D\in {\mathcal O}_X(D)$ be a holomorphic section of $[D]$ defining $D$ and by rescaling $h_D$ if necessary we assume $\sup_X |s_D|^2_{h_D}\le 1$.

\smallskip

To ease the notations, throughout this section we  will denote by $\tilde g$ and $g$ (omitting the subscript $t$ in $g_t$) the associated metrics of $\omega_X$ and $\omega_t$, respectively. 

 We will omit the subscript $t$ in $\varphi_t$ which solves \neweqref{eqn:MA s3} and simply write it as  $\varphi$. The function $F_{\omega_t}$ will be simply written as $F$, since $\omega_t$ is a fixed metric.

Recall that we denote ${\mathcal G}_t$ the positive Green's function in \neweqref{eqn:positive Green}, which differs from $G_t$ by a uniform constant.

\subsection{Gradient estimates} 

We will prove Theorem \ref{thm:gradient} in this subsection. Fix a constant $p>n$.

\smallskip

The lemma  below follows from straightforward calculations, so we omit the proof.
\begin{lemma}\label{lemma gradient trivial}
Suppose $\varphi$ satisfies \neweqref{eqn:MA s3}. We have 
\begin{equation}\label{eqn:grad 3.1}
\Delta_{g} |\nabla \varphi|_{\tilde g}^2\ge 2 Re\langle \nabla F, \bar \nabla \varphi \rangle_{\tilde g} + g^{i\bar j} \tilde g^{k\bar l} ( \varphi_{ki} \varphi_{\bar j \bar l} + \varphi_{k\bar j} \varphi_{i\bar l}  ) - 2 K \tr_{g} \tilde g |\nabla \varphi|^2_{\tilde g},
\end{equation}
where $-K$ is a lower bound of the bisectional curvature of the fixed metric $\omega_X$ and $\varphi_{ki}$ denote the second-order covariant derivatives of $\varphi$ with respect to $\omega_X$.
\end{lemma}
\begin{lemma}\label{lemma 3.2}
The following inequality holds on $X$:
\begin{equation}\label{eqn:grad 3.2}
\Delta_g H\ge 2 e^{-\lambda \hat \varphi} Re\langle \nabla F, \bar \nabla \varphi\rangle_{\tilde g} + H \tr_{g} \tilde g- (2+n) \lambda H - C H^{1/2} - C_s H^{1/2 } \tr_g \tilde g.
\end{equation}
where $\hat \varphi = \varphi - \varepsilon_0 \log |s_D|_{h_D}^2$, and $H = e^{-\lambda \hat\varphi} |\nabla \varphi|^2_{\tilde g}$ for suitable $\lambda>0$ depending only on $\varepsilon_0, \delta_0$ and $\omega_X$. The constants $C>0$ and $C_s>0$ in \neweqref{eqn:grad 3.2} are both uniform.
\end{lemma}
\noindent{\em Proof. }
As above, we denote by $H = e^{-\lambda \hat\varphi} |\nabla \varphi|^2_{\tilde g}$ for some $\lambda>0$ to be determined later. We calculate on $X\backslash D$
\bea
\Delta_g H \nonumber &= &\nonumber e^{-\lambda \hat\varphi} \Delta_g |\nabla \varphi|_{\tilde g}^2 + 2 e^{-\lambda \hat \varphi} Re\langle -\lambda \nabla \hat \varphi, \bar \nabla |\nabla \varphi|_{\tilde g}^2 \rangle_{g} + H ( -\lambda \Delta_g \hat \varphi + \lambda^2 |\nabla \hat \varphi|_g  )\\
&\ge \label{eqn:3.5} &  e^{-\lambda \hat\varphi} \Big(2 Re\langle \nabla F, \bar \nabla \varphi \rangle_{\tilde g} + g^{i\bar j} \tilde g^{k\bar l} ( \varphi_{ki} \varphi_{\bar j \bar l} + \varphi_{k\bar j} \varphi_{i\bar l}  ) - 2 K \tr_{g} \tilde g |\nabla \varphi|^2_{\tilde g}  \Big) \\
& & \nonumber + 2 e^{-\lambda \hat \varphi} Re\langle -\lambda \nabla \hat \varphi, \bar \nabla |\nabla \varphi|_{\tilde g}^2 \rangle_{g} + H ( -\lambda \Delta_g \hat \varphi + \lambda^2 |\nabla \hat \varphi|_g  ).
\eea
We now perform the calculation at a fixed point $x_0\in X\backslash D$ and choose normal coordinates at $x_0$ relative to $\tilde g$, such that $\tilde g_{i\bar j}(x_0) = \delta_{ij}$, $d\tilde g_{i\bar j}(x_0) = 0$ and $g_{i\bar j}(x_0) = g_{i\bar i} \delta_{ij}$ is diagonal. Then at $x_0$ we have
\bea
\nonumber && 2  e^{-\lambda \hat \varphi} Re\langle -\lambda \nabla \hat \varphi, \bar \nabla |\nabla \varphi|_{\tilde g}^2 \rangle_{g}  =  -2 \lambda e^{-\lambda \hat \varphi} Re( g^{i\bar i}( \hat \varphi_i \varphi_{k\bar i} \varphi_{\bar k}  + \hat \varphi_i \varphi_{\bar k\bar i} \varphi_k   ) )\\
\nonumber &\ge & -2 \lambda e^{-\lambda\hat\varphi} Re( g^{i\bar i} \hat \varphi_i\varphi_{i\bar i} \varphi_{\bar i}  ) - \lambda^2 e^{-\lambda \hat \varphi} |\nabla \hat \varphi|_{g}^2 |\nabla \varphi|_{\tilde g}^2  - e^{-\lambda \hat \varphi} g^{i\bar i} \varphi_{\bar k \bar i} \varphi_{k i } \\
\nonumber &= & -2 \lambda e^{-\lambda\hat\varphi} Re( g^{i\bar i} \hat \varphi_i(g_{i\bar i} - 1) \varphi_{\bar i}  ) - \lambda^2 H |\nabla \hat \varphi|_{g}^2  - e^{-\lambda \hat \varphi} g^{i\bar i} \varphi_{\bar k \bar i} \varphi_{k i } \\
\nonumber &= & -2 \lambda e^{-\lambda\hat\varphi} \langle \nabla \hat \varphi, \nabla \varphi\rangle_{\tilde g} + 2 \lambda e^{-\lambda \hat\varphi} \langle \nabla \hat \varphi, \nabla \varphi\rangle_g   - \lambda^2 H |\nabla \hat \varphi|_{g}^2  - e^{-\lambda \hat \varphi} g^{i\bar i} \varphi_{\bar k \bar i} \varphi_{k i } .
\eea
And  by \neweqref{eqn:Kodaira} we have
$$- \Delta_g \hat\varphi = \tr_{g} (\chi + t\omega_X - \omega_t - \varepsilon_0 \ric(h_D))\ge \delta_0 \tr_{g} \tilde g - n$$
Substituting the above two (in)equalities  to \neweqref{eqn:3.5}, we obtain that at $x_0$
\bea
\Delta_g H  \nonumber &\ge &  e^{-\lambda \hat\varphi} \Big(2 Re\langle \nabla F, \bar \nabla \varphi \rangle_{\tilde g} + g^{i\bar j} \tilde g^{k\bar l}  \varphi_{k\bar j} \varphi_{i\bar l} - 2 K \tr_{g} \tilde g |\nabla \varphi|^2_{\tilde g}  \Big)  \nonumber + \lambda H ( \delta_0\tr_g \tilde g - n )\\
\nonumber &&  -2 \lambda e^{-\lambda\hat\varphi} \langle \nabla \hat \varphi, \nabla \varphi\rangle_{\tilde g} + 2 \lambda e^{-\lambda \hat\varphi} \langle \nabla \hat \varphi, \nabla \varphi\rangle_g \nonumber\\
&\ge & \nonumber 2 e^{-\lambda \hat \varphi} Re \langle \nabla F, \bar \nabla \varphi \rangle_{\tilde g} + (\lambda \delta_0 - 2K) (\tr_g \tilde g ) H - \lambda n H - 2 \lambda e^{-\lambda \hat\varphi} |\nabla \varphi|_{\tilde g}^2 + 2 \lambda e^{-\lambda\hat\varphi} |\nabla \varphi|_g^2 \\
& & \label{eqn:3.6} - 2\varepsilon_0 \lambda e^{-\lambda \varphi} |s_D|_{h_D}^{2 \lambda \varepsilon_0 -1 } |\nabla_{h_D} s_D|_{\tilde g, h_D} |\nabla \varphi|_{\tilde g} - 2 \varepsilon_0 \lambda  e^{-\lambda \varphi} |s_D|_{h_D}^{2 \lambda \varepsilon_0 -1 } |\nabla_{h_D} s_D|_{ g, h_D} |\nabla \varphi|_{ g},
\eea
where we denote by $\nabla_{h_D} $ the Chern connection of the hermitian metric $h_D$ on the line bundle $[D]$. We choose $\lambda>1$ large enough so that $\lambda \delta_0 > 2K+10$ and $2 \lambda \varepsilon_0>10$. Then the first term in \neweqref{eqn:3.6} satisfies
\begin{equation}\label{eqn:3.7}
- 2\varepsilon_0 \lambda e^{-\lambda \varphi} |s_D|_{h_D}^{2 \lambda \varepsilon_0 -1 } |\nabla_{h_D} s_D|_{\tilde g, h_D} |\nabla \varphi|_{\tilde g}\ge - C e^{-\lambda \varphi/2} |s_D|_{h_D}^{\lambda \varepsilon_0 - 1} H^{1/2}\ge -C H^{1/2},
\end{equation}
because $|\nabla _{h_D} s_D|_{\tilde g, h_D}\le C$ and $\| \varphi\|_{L^\infty}\le C$.  Similarly the last term in \neweqref{eqn:3.6} satisfies
\begin{equation}\label{eqn:3.8}
 - 2 \varepsilon_0 \lambda  e^{-\lambda \varphi} |s_D|_{h_D}^{2 \lambda \varepsilon_0 -1 } |\nabla_{h_D} s_D|_{ g, h_D} |\nabla \varphi|_{ g} \ge - C_s H^{1/2} \tr_g \tilde g,
\end{equation}
for some uniform constant $C_s>0$. Plugging \neweqref{eqn:3.7} and \neweqref{eqn:3.8} into \neweqref{eqn:3.6}, we obtain that at $x_0\in X\backslash D$
\begin{equation}\label{eqn:3.91}
\Delta_g H\ge 2 e^{-\lambda \hat \varphi} Re\langle \nabla F, \bar \nabla \varphi\rangle_{\tilde g} + H \tr_{g} \tilde g- (2+n) \lambda H - C H^{1/2} - C_s H^{1/2 } \tr_g \tilde g.
\end{equation}
Since $x_0\in X\backslash D$ is arbitrary, $X\backslash D\subset X$ is clearly dense and both sides of \neweqref{eqn:3.91} are smooth, we see that \neweqref{eqn:3.91} holds globally on $X$. 

\medskip

\noindent{\em Proof of Theorem \ref{thm:gradient}.} Let $C_s>0$ be the constant in Lemma \ref{lemma 3.2} and denote $\Lambda = C_s^2+1$. From Lemma \ref{lemma 3.2}, we have for some uniform constant $C>0$
\begin{equation}\label{eqn:3.9}
\Delta_g H\ge -C |\nabla F|_{\tilde g} H^{1/2} + (H^{1/2} - C_s) H^{1/2} \tr_{g} \tilde g- C H - C 
\end{equation}
We consider the convex and monotonically increasing  function $\Phi_\delta(x) = \frac{1}{2}( \sqrt{x^2 + \delta} + x  ) + \Lambda$ for $\delta>0$. It is clear that
$$\Phi_\delta(H-\Lambda) \to \max(H, \Lambda),\mbox{ as } \delta\to 0.$$
We also have $0\le \Phi_\delta'(x)\le 1$, $\Phi_\delta''(x)\ge 0$ for all $x\in{\mathbb R}$ and
\begin{equation}\label{eqn:3.10}0\le \Phi_\delta'(x) = \frac{1}{2} \frac{\delta}{\sqrt{x^2 + \delta} (\sqrt{x^2 + \delta} - x)} \le \frac{\delta}{2},\mbox{ when }x\le -1.
\end{equation}
We denote $\hat H = H - \Lambda$ and calculate
\bea
\Delta_g \Phi_\delta(\hat H) &= \nonumber &\, \Phi_\delta' \Delta_g H + \Phi_\delta'' |\nabla H|_g^2\\
& \ge &\, \nonumber  \Phi_\delta'(\hat H) \Big( -C |\nabla F|_{\tilde g} H^{1/2} + (H^{1/2} - C_s) H^{1/2} \tr_{g} \tilde g- C H - C\Big) \chi_{\{\hat H \ge -1\}} \\
& &\, \nonumber +  \Delta_g H \cdot \Phi_\delta'(\hat H) \chi_{\{\hat H<-1\}} \\
& \ge & \, \nonumber \Phi_\delta'(\hat H) \Big( -C |\nabla F|_{\tilde g} H^{1/2} - C H - C\Big) \chi_{\{\hat H \ge -1\}}  -\delta  |\Delta_g H|  \chi_{\{\hat H<-1\}} 
\eea
where the second inequality follows since on $\{\hat H \ge -1\}$, $H\ge \Lambda - 1 \ge  C_s^2$, and on $\{\hat H<-1\}$, we have $\Phi_\delta'(\hat H)\le \delta$ by \neweqref{eqn:3.10}. Applying Green's formula to the {\em smooth} function $\Phi_\delta(\hat H)$ we obtain for any $x\in X$
\bea
\Phi_\delta(\hat H)(x)& = & \label{eqn:3.12} \, \frac{1}{V_t}\int_X \Phi_\delta(\hat H) \omega_t^n + \int_X {\mathcal G}_t(x,\cdot)(-  \Delta_{g} \Phi_\delta(\hat H )) \omega_t^n\\
 &\le & \nonumber \, \frac{1}{V_t}\int_X \Phi_\delta(\hat H) \omega_t^n + \delta \int_{\{\hat H<-1\}} {\mathcal G}_t(x,\cdot)  |\Delta_g H| \omega_t^n  \\
& & \, \nonumber+ \int_{\{\hat H\ge -1\}} {\mathcal G}_t(x,\cdot)  \Phi_\delta'(\hat H) \Big(C |\nabla F|_{\tilde g} H^{1/2} + C H + C\Big) \omega_t^n.
\eea
Letting $\delta\to 0$ in \neweqref{eqn:3.12}, we get (denoting $H_\Lambda = \max\{H, \Lambda\}$)
\bea\label{eqn:3.13}
H_\Lambda(x) \le \frac{1}{V_t}\int_X H_\Lambda \omega_t^n +  \int_{X} {\mathcal G}_t(x,\cdot) \Big(C |\nabla F|_{\tilde g} H^{1/2} + C H + C\Big) \omega_t^n.
\eea
Let $x_1\in X$ be a maximum point of $H_\Lambda$. Then we obtain from \neweqref{eqn:3.13} that
{\small
\bea
H_{\Lambda,\max} &= \nonumber & \, H_\Lambda (x_1)
\le \frac{ c_t H_{\Lambda, \max}^{1-\eta}}{V_t} \int_X H_\Lambda^\eta e^F \omega_X^n \\
&& \, \nonumber  + C H_{\Lambda,\max}^{1/2} \int_X {\mathcal G}_t(x_1, \cdot) |\nabla F|_{\tilde g} \omega_t^n + C H_{\Lambda,\max}^{1-\eta} \int_X {\mathcal G}_t(x_1, \cdot) H^{\eta} \omega_t^n + C\\
& \le & \, \label{eqn:3.14}  \frac{ c_t H_{\Lambda, \max}^{1-\eta}}{V_t} \Big( \int_X  e^{(1+\epsilon)F} \omega_X^n\Big)^{\frac{1}{1+\epsilon}} \Big( \int_X H_\Lambda^{\eta (1+\epsilon)/\epsilon} \omega_X^n \Big)^{\frac{\epsilon}{1+\epsilon}} +C  \\
& & \, \nonumber  + C H_{\Lambda,\max}^{1/2}\Big( \int_X {\mathcal G}_t(x_1, \cdot)^{p^*} \omega_t^n \Big)^{{1/}{p^*}} \Big( \int_X |\nabla F|_{\tilde g}^p c_t e^F \omega_X^n    \Big)^{{1/}{p}}\\
& & \, \nonumber + C H_{\Lambda,\max}^{1-\eta}  \Big( \int_X {\mathcal G}_t(x_1, \cdot)^{p_0} \omega_t^n \Big)^{1/p_0}\Big( \int_X H^{\eta p_0^*} \omega_t^n   \Big)^{1/p_0^*},
\eea
}%
where we fix a number $p_0\in (1, \frac{n}{n-1})$. By assumption $p>n$ so we have $p^*<\frac{n}{n-1}$. If furthermore we choose $\eta>0$ small (depending on only $\epsilon$ and $p_0$) such that $\eta (1+\epsilon)/\epsilon \le 1$ and $\eta p_0^*\le \epsilon/(1+\epsilon)$, then by Theorem \ref{thm:main2}, and Lemma \ref{lemma L2} below it follows that all the integrals involved in \neweqref{eqn:3.14} are bounded uniformly from above. It then follows that 
$$H_{\Lambda,\max}\le C H_{\Lambda,\max}^{1-\eta} + CH_{\Lambda,\max}^{1/2} + C.$$
By Young's inequality we immediately derive the uniform upper bound of $H_{\Lambda, \max}$. Hence the upper bound of $H = e^{-\lambda \varphi}|s_D|^{2\varepsilon_0\lambda}_{h_D} |\nabla \varphi|^2_{\omega_X}$. The proof of  Theorem \ref{thm:gradient} is completed.

\medskip

It only remains to show the integrals of $H$ are bounded.

\begin{lemma}
\label{lemma L2}
There is a uniform constant $C>0$ such that 
$$\int_X H \omega_X^n = \int_X e^{-\lambda \varphi} |s_D|^{2\lambda \varepsilon_0}_{h_D} |\nabla \varphi|^2_{\omega_X} \omega_X^n \le C,$$
and when integrated against $\omega_t^n$ we have
$$\int_X H^{\frac{\epsilon}{1+\epsilon}} \omega_t^n \le C.$$
\end{lemma}
\noindent{\em Proof. }Recall we write $\hat \varphi = \varphi - \varepsilon_0 \log |s_D|^2_{h_D}$. 
For any small $\delta>0$, we denote the super-level set of $|s_D|_{h_D}$,  $E_\delta = \{|s_D|_{h_D} \ge \delta\}$. We write $\hat \omega = \chi - \varepsilon_0\ric(h_D)\ge \delta_0 \omega_X$ as in \neweqref{eqn:Kodaira}. Observe that on $E_\delta$ the following equation holds
\bea
   \omega_t^n - (\hat \omega + t \omega_X)^n &=&  (\hat \omega + t \omega_X + \ddbar \hat \varphi)^n - (\hat \omega + t \omega_X)^n \nonumber\\
&  = & \, \nonumber  \ddbar \hat \varphi  \wedge ( \omega_t^{n-1} + \cdots + (\hat \omega + t \omega_X)^{n-1}   ) . 
\eea
Multiplying both sides by $e^{-\lambda \hat \varphi}$ 
 and integrating over $E_\delta$, we obtain by integration by parts
 \bea
&& \, \frac{\lambda}{n}\int_{E_\delta} e^{-\lambda \hat \varphi} |\nabla \hat\varphi|^2_{\hat \omega} \hat \omega^n\le    \nonumber {\lambda}  \int_{E_\delta} e^{-\lambda \hat \varphi} i\partial \hat \varphi \wedge \bar \partial \hat \varphi  \wedge ( \omega_t^{n-1} + \cdots + (\hat \omega + t \omega_X)^{n-1}   ) \\
&= & \,  \nonumber \int_{\partial E_\delta} e^{-\lambda \hat\varphi} i\bar \partial \hat\varphi  \wedge ( \omega_t^{n-1} + \cdots + (\hat \omega + t \omega_X)^{n-1}   ) + \int_{E_\delta} e^{-\lambda \hat\varphi} (\omega_t^n - (\hat \omega + t\omega_X)^n)\\
&\le & \,  \nonumber \int_{\partial E_\delta} e^{-\lambda \hat\varphi} i\bar \partial \hat\varphi  \wedge ( \omega_t^{n-1} + \cdots + (\hat \omega + t \omega_X)^{n-1}   ) + C,
 \eea
 since the function $e^{-\lambda \hat\varphi} = e^{-\lambda\varphi} |s_D|^{2\lambda\varepsilon_0}_{h_D}$ is uniformly bounded. 
\vskip .1in
We rewrite the boundary integral as follows.
 \bea
&& \nonumber\, \int_{\partial E_\delta} e^{-\lambda \hat\varphi} i\bar \partial (\varphi - \varepsilon_0 \log |s_D|^2_{h_D}) \wedge ( \omega_t^{n-1} + \cdots + (\hat \omega + t \omega_X)^{n-1}   )\\
& = & \     \nonumber 
\int_{\partial E_\delta} e^{-\lambda \varphi}|s_D|^{2\lambda\varepsilon_0-2}_{h_D}\,  \eta  \wedge ( \omega_t^{n-1} + \cdots + (\hat \omega + t \omega_X)^{n-1}   ) \\
&= & \,  \nonumber \delta^{2\lambda\epsilon_0-2}\int_{E_\delta}\, 
d(e^{-\lambda \varphi}\,  \eta)  \wedge ( \omega_t^{n-1} + \cdots + (\hat \omega + t \omega_X)^{n-1}   ) 
 \eea
Here $\eta=i|s_D|^2_{_{h_D}}\bar\partial (\varphi - \varepsilon_0 \log |s_D|^2_{h_D})$ is a smooth $1$-form on $X$ and in the last line we have integrated by parts again. Now the last integrand is a smooth $(n,n)$ form on $X$ which is independent of $\delta$ and thus the integral remains bounded as $\delta\rightarrow 0$. On the other hand, $2\lambda\epsilon_0-2\geq 10-2>0$ so the last line tends to zero as $\delta\rightarrow 0$.
 \vskip .1in
 Letting $\delta\to 0$ we get from the equivalence of the fixed metrics $\omega_X$ and $\hat \omega$ that
 $$\int_{X} e^{-\lambda \hat\varphi}|\nabla \hat\varphi|_{\omega_X}^2 \omega_X^n\le C. $$
 The first inequality in the lemma then follows from the following triangle inequality
 $$\int_{X} e^{-\lambda \hat\varphi}|\nabla \varphi|_{\omega_X}^2 \omega_X^n\le 2\int_{X} e^{-\lambda \hat\varphi}|\nabla \hat\varphi|_{\omega_X}^2 \omega_X^n   + 2 \varepsilon_0^2\int e^{-\lambda \varphi} |s_D|^{2\varepsilon_0 \lambda }_{h_D} |\nabla \log |s_D|_{h_D}^2|^2_{\omega_X} \omega_X^n\le C.     $$ The second inequality follows from H\"older inequality 
 $$\int_X H^{\frac{\epsilon}{1+\epsilon}} \omega_t^n  \le c_t \Big(\int_X e^{(1+\epsilon) F} \omega_X^n \Big)^{1/(1+\epsilon)} \Big(\int_X H \omega_X^n \Big)^{\epsilon/(1+\epsilon)}\le C.$$%
 The proof of Lemma \ref{lemma L2} is completed.
 
 \medskip

We recall the gradient estimate proved in Theorem \ref{thm:gradient}  for complex Monge-Amp\`ere equations with a fixed background metric.
\begin{equation}\label{eqn:fixed MA}
(\omega_X + \ddbar \varphi)^n = e^F \omega_X^n,\quad \sup_X \varphi=0.
\end{equation}
The conditions on the sets ${\mathcal M}'(N,\epsilon, \gamma)$, ${\mathcal M}''(N,\epsilon, \gamma)$ or $\tilde {\mathcal M}''(N,\epsilon, \gamma)$ state that
\begin{equation}\label{eqn:F condition}
\| e^F\|_{L^{1+\epsilon}(X,\omega_X^n)}\le N,\quad \mbox{and}
\end{equation}
{\small
\begin{equation}\label{eqn:F condition 1}
\sup e^{-F}\le \gamma,\, \mbox{or }\int_X (e^{-F} + |\Delta_{\omega_X} e^{-F}|) \omega_X^n \le \gamma, \,\mbox{or }\int_X (e^{-F} + |\nabla e^{-F}|_{\omega_X}^2) \omega_X^n \le \gamma
\end{equation}
}
A corollary of Theorem \ref{thm:gradient} on equation \neweqref{eqn:fixed MA} states that 

\begin{corollary}
\label{corollary for gradient} Fix a constant $p>n$.
Let $\varphi$ solve the equation \neweqref{eqn:fixed MA}. Suppose $F$ satisfies \neweqref{eqn:F condition} and \neweqref{eqn:F condition 1}, then the following gradient estimate of $\varphi$ w.r.t. $\omega_X$ holds
\begin{equation}\label{eqn:gradient fixed}
\sup_X |\nabla \varphi|^2_{\omega_X} \le C,
\end{equation}
where $C>0$ depends on $n, p , \omega_X, N, \epsilon, \gamma$, and  additionally $\int_X |\nabla F|_{\omega_X}^p e^F \omega_X^n$.
\end{corollary}
Corollary \ref{corollary for gradient} follows from Theorem \ref{thm:gradient} as a particular case if we choose $\chi = \omega_X/2$, $t=1/2$, $D = {\mathbf{0}}$ (the trivial divisor), $s_D\equiv 1$ and $h_D\equiv 1$.

\smallskip

 We remark that besides the integral bounds on $e^{\pm F}$ and $|\nabla F|$, the gradient estimate of $\varphi$ can be made to be  independent of the pointwise bounds $\sup_X e^F$ and $\inf_X e^F$. The exponent $p>n$ in the $L^p(X, e^F \omega_X^n)$-bound of $|\nabla F|_{\omega_X}$ is also {\em sharp} in the sense that the gradient estimate may fail if $p<n$ and all other conditions on $F$ are still valid, as the following example shows. We do not know whether the statement holds or not when $p=n$.
 
 \medskip

\noindent{\bf Example 3.1.} Let $0\in {\mathbb C}^n\subset {\mathbb {CP}}^n$ and $z=(z_1,\ldots, z_n)$ be the natural coordinates on ${\mathbb C}^n$. Fix a number $a\in (0,1)$ and a small $\delta\in (0,1/100)$. We consider the function
$$
\varphi_\delta = \left\{\begin{array}{lll}
&(|z|^2+ \delta)^{a}, & \mbox{ if }|z|<\zeta \\
&\widetilde{\max}\Big\{ (|z|^2+ \delta)^{a}, 2 \log (1+ |z|^2)  \Big\}, & \mbox{ if } \zeta\le |z| \le 1\\
& 2\log (1+|z|^2), & \mbox{ if } |z|\ge 1,
\end{array} \right.
$$
where $\zeta>0$ is a constant (independent of $\delta$) such that $(|z|^2+ \delta)^{a} = 2 \log (1+ |z|^2)$ for some $|z|\in (\zeta, 1)$. The metric $\omega_\delta = \ddbar \varphi_\delta$ on ${\mathbb C}^n$ can be naturally extended to a smooth K\"ahler metric on $X = {\mathbb{CP}}^n$. We express the metric $\omega_\delta$ locally near $0$
$$
\omega_\delta =   \, \sum_{i,j} a \frac{\delta_{ij} + (a-1){\bar z_i z_j/}{(|z|^2 + \delta)}}{ (|z|^2 + \delta)^{1-a}  } \sqrt{-1} dz_i\wedge d\bar z_j
$$
so  near $0$ we have 
$$\omega_\delta^n = a^n \frac{(a|z|^2 + \delta)/ (|z|^2 + \delta)}{(|z|^2 + \delta)^{n -an}} \omega_{\mathbb C^n}^n = :e^F \omega_{\mathbb C^n}^n.$$
If $a<1/2$, we can choose $0<p< 2an <n$. By straightforward calculations, we see that  near $0$
\begin{equation}\label{eqn:example}
|\nabla F|^p_{\omega_{{\mathbb C}^n}} e^F \le C \frac{|z|^p}{|z|^{2p + 2n - 2an}} = C |z|^{-p + 2an - 2n}
\end{equation}
and the function on the right-hand side is integrable near $z=0$. Since $\omega_{{\mathbb C}^n}$ is equivalent to $\omega_{FS}$ (the Fubini-Study metric on ${\mathbb {CP}}^n$) near $0$, it follows easily that $\int_X |\nabla F|_{\omega_{FS}}^p e^F \omega_{FS}^n$ is uniformly bounded above (independent of $\delta>0$). Moreover, the other conditions \neweqref{eqn:F condition} and \neweqref{eqn:F condition 1} of $F$ are satisfied for certain $\epsilon, N$ and $\gamma$. Here we have viewed $F$ as a smooth function on the whole manifold ${\mathbb {CP}}^n$. However, near $0$ (e.g. at $|z|^2 = \delta$)
$$|\nabla \varphi_\delta|_{\omega_{\mathbb C^n}}\sim \frac{|z|}{(|z|^2 + \delta)^{1-a}} \sim \delta^{a - \frac{1}{2}}$$
blows up as $\delta\to 0$. Therefore when $p<n$, the integral $\int_X |\nabla F|_{\omega_X}^p e^F \omega_X^n$ is not enough to conclude the gradient estimate of $\varphi$ which satisfies  \neweqref{eqn:fixed MA}.

\subsection{$C^2$ estimate} 

We consider the $C^2$ estimates in this subsection and give the proof of Theorem \ref{thm:C2}. %
We continue to use the same notations as in the previous subsection. Fix a number $p>2n$. 

\smallskip

Let $\varphi_t$ be the solution to the equation \neweqref{eqn:MA s3}. We again omit the subscript $t$ in $\varphi_t$.

\newcommand{\innpro}[1]{\langle#1\rangle}

\begin{lemma}\label{lemma Yau C2}
The following holds for any $t\in (0,1]$:
\bea
\Delta_g \tr_{\tilde g} g&\ge&\, \label{eqn:YauC2} \Delta_{\tilde g} F - K \tr_g \tilde g \cdot \tr_{\tilde g} g + \tilde g^{k\bar l} \nabla^{\tilde g}_{ k} g_{i\bar q} \bar \nabla^{\tilde g}_{ \bar l} g_{p\bar j} g^{p \bar q} g^{i\bar j} \\
&\ge &\, \nonumber \Delta_{\tilde g} F - K \tr_g \tilde g \cdot \tr_{\tilde g} g + \frac{|\nabla \tr_{\tilde g} g|^2_{g}}{\tr_{\tilde g} g},
\eea
where as before $-K$ is a lower bound of the bisectional curvature of $\tilde g = \omega_X$, and $\nabla^{\tilde g}$ denotes the covariant derivatives with respect to $\tilde g$.
\end{lemma}
Lemma \ref{lemma Yau C2} follows from standard calculations in \cite{Y}, so we omit the proof.
\begin{lemma}\label{lemma C2}
For some $\mu>0$ depending only on $\chi$, and $ \omega_X$, we have
$$\Delta_g Q \ge e^{-\mu\hat \varphi } \Delta_{\tilde g} F - n\mu Q,$$
where we have written $Q = e^{-\mu \hat\varphi} \tr_{\tilde g} g$ and $\hat \varphi = \varphi - \varepsilon_0 \log |s_D|^2_{h_D}$.
\end{lemma}
\noindent{\em Proof. }
We fix a constant $\mu>0$ to be determined later. We calculate using Lemma \ref{lemma Yau C2}
\bea
\Delta_g Q &= & \, \nonumber e^{-\mu \hat\varphi} \Delta_g \tr_{\tilde g} g - 2\mu e^{-\mu \hat\varphi} Re \innpro{\nabla\hat \varphi, \bar \nabla \tr_{\tilde g} g }_{g} + Q(- \mu \Delta_g \hat \varphi + \mu^2 |\nabla \hat\varphi|_g^2)\\
&\ge & \, \label{eqn:3.18} e^{-\mu \hat\varphi}\big( \Delta_{\tilde g} F - K \tr_g \tilde g \cdot \tr_{\tilde g} g + \frac{|\nabla \tr_{\tilde g} g|^2_{g}}{\tr_{\tilde g} g} \big) - 2\mu e^{-\mu \hat\varphi} Re \innpro{\nabla\hat \varphi, \bar \nabla \tr_{\tilde g} g }_{g}  \\
&&\,\, \nonumber + \mu Q\Big[- n + t \cdot  \tr_g \tilde g +  \tr_{g}(\chi - \varepsilon_0 \ric(h_D)) \Big] + \mu^2 e^{-\mu \hat\varphi} \tr_{\tilde g} g \cdot |\nabla \hat\varphi|_g^2.
\eea
Applying Cauchy-Schwarz inequality the following holds
$$- 2\mu e^{-\mu \hat\varphi} Re \innpro{\nabla\hat \varphi, \bar \nabla \tr_{\tilde g} g }_{g} \ge  - \mu^2 e^{-\mu \hat \varphi} \tr_{\tilde g} g |\nabla \hat \varphi|_g^2 - e^{-\mu\hat\varphi}\frac{|\nabla \tr_{\tilde g} g|^2_g}{\tr_{\tilde g} g},$$
plugging this into \neweqref{eqn:3.18} and applying \neweqref{eqn:Kodaira} we obtain 
\bea
\Delta_g Q& \ge& \nonumber \, e^{-\mu\hat \varphi } \Delta_{\tilde g} F - n\mu Q + (\mu \delta_0 - K) e^{-\mu\hat \varphi} \tr_g \tilde g\cdot \tr_{\tilde g} g\\
 &\ge& \, \nonumber e^{-\mu\hat \varphi } \Delta_{\tilde g} F - n\mu Q,
\eea
if we choose $\mu>0$ such that $\mu\delta_0 - K \ge 1$.

\medskip

\noindent{\em Proof of Theorem \ref{thm:C2}.}  By Lemma \ref{lemma C2} and Green's formula we have for any $x\in X$
\bea
Q(x) &= & \nonumber \; \frac{1}{V_t} \int_X Q \omega_t^n + \int_X {\mathcal G}_t(x,\cdot) (-\Delta_g Q) \omega_t^n\\
&\le &  \label{eqn:3.19} \; \frac{1}{V_t} \int_X Q \omega_t^n + \int_X {\mathcal G}_t(x,\cdot) (- e^{-\mu \hat \varphi} \Delta_{\tilde g} F  + n \mu Q) \omega_t^n.
\eea
Let $x_0$ be a maximum point of $Q$. We apply \neweqref{eqn:3.19} at the point $x_0$, then
{\small
\bea
Q_{\max} &= & \; \nonumber Q(x_0) \le  \frac{1}{V_t} \int_X Q \omega_t^n + \int_X {\mathcal G}_t(x_0,\cdot) (- e^{-\mu \hat \varphi} \Delta_{\tilde g} F  + n \mu Q) \omega_t^n\\
&\le & \; \label{eqn:3.20}\frac{Q_{\max}^{1-\eta}}{V_t} \int_X Q^\eta \omega_t^n + n \mu Q_{\max}^{1-\eta} \int_X {\mathcal G}_t(x_0,\cdot) Q^\eta\omega_t^n  + c_t \int_X {\mathcal G}_t(x_0,\cdot) (- e^{-\mu \hat \varphi} \Delta_{\tilde g} F)e^F \omega_X^n.
\eea
}
We observe that by H\"older inequality
{\small
\begin{equation}\label{eqn:bound Q}
\int_X Q^{\frac{\epsilon}{1+\epsilon}} \omega_t^n =c_t \int_X Q^{\frac{\epsilon}{1+\epsilon}}e^F \omega_X^n  \le c_t\| e^F\|_{L^{1+\epsilon}} \Big( \int_X e^{-\mu \hat \varphi } (\tr_{\omega_X} \omega_t) \omega_X^n\Big)^{(1+\epsilon)/\epsilon}\le C,
\end{equation}
}%
since $e^{-\mu\hat \varphi} = e^{-\mu \varphi} |s_D|_{h_D}^{2\mu \varepsilon_0}\le C$.  By H\"older inequality and \neweqref{eqn:bound Q}, the integral in the first term in \neweqref{eqn:3.20} is bounded if $0<\eta\le \frac{\epsilon}{1+\epsilon}$. The integral in the second term in \neweqref{eqn:3.20} satisfies
$$\int_X {\mathcal G}_t(x_0,\cdot) Q^\eta\omega_t^n \le \Big(\int_X {\mathcal G}_t(x_0, \cdot)^{p_0} \omega_t^n \Big)^{1/p_0}  \Big(\int_X Q^{p_0^* \eta} \omega_t^n \Big)^{1/p_0^*}\le C$$
if we choose  $p_0\in (1, \frac{n}{n-1})$ and $\eta>0$ small so that $p_0^* \eta\le \frac{\epsilon}{1+\epsilon}$. The inequality above then follows from \neweqref{eqn:bound Q} and Theorem \ref{thm:main2}.

\newcommand{\xk}[1]{\big(#1\big)}
\smallskip

To deal with the last integral in \neweqref{eqn:3.20}, we apply integration by parts to obtain 
{\small
\bea
&& \; \nonumber c_t\int_X {\mathcal G}_t(x_0,\cdot) (- e^{-\mu \hat \varphi} \Delta_{\tilde g} F)e^F \omega_X^n\\
&= & \; \nonumber c_t\int_{X} \Big( {\mathcal G}_t(x_0,\cdot) e^{-\mu \hat \varphi} |\nabla F|_{\tilde g}^2  + {\mathcal G}_t(x_0,\cdot) \innpro{\nabla F, \bar \nabla e^{-\mu \hat\varphi}} _{\tilde g}  + e^{-\mu\hat\varphi} \innpro{\nabla {\mathcal G}_t(x_0,\cdot), \nabla F}_{\tilde g} \Big)e^F \omega_X^n\\
\label{eqn:3.21}&\le & \; c_t \int_X\Big (2 {\mathcal G}_t(x_0,\cdot) e^{-\mu \hat \varphi} |\nabla F|^2_{\tilde g}  + \mu^2 {\mathcal G}_t(x_0,\cdot) e^{-\mu\hat\varphi} |\nabla \hat\varphi |_{\tilde g}^2 + e^{-\mu\hat\varphi} |\nabla F|_{\tilde g} |\nabla {\mathcal G}_t(x_0,\cdot)|_{\tilde g}\Big) e^F.
\eea 
}
The first integral in \neweqref{eqn:3.21} satisfies (recall $p>2n$)
$$ 2 c_t\int_X {\mathcal G}_t(x_0,\cdot) e^{-\mu \hat \varphi} |\nabla F|^2_{\tilde g} e^F \le C \Big( \int_X {\mathcal G}_t(x_0,\cdot)^{(p/2)^*} \omega_t^n \Big)^{1/(p/2)^*} \cdot \Big( \int_X |\nabla F|_{\tilde g}^{p} e^F \omega_X^n   \Big)^{2/p}$$
which is bounded above uniformly since $(p/2)^* = \frac{p/2}{p/2-1}< \frac{n}{n-1}$ we can apply Theorem \ref{thm:main2} and the assumption on $F$ in Theorem \ref{thm:C2}. The second integral in \neweqref{eqn:3.21} is also bounded because of Theorem \ref{thm:main2} and $e^{-\mu\hat\varphi} |\nabla \hat \varphi|_{\tilde g}^2\le C$ which follows from Theorem \ref{thm:gradient} (we may assume $\mu> A$, where $A>0$ is the constant in Theorem \ref{thm:gradient}). We now deal with the last integral in \neweqref{eqn:3.21}.  We calculate
{\small
\bea
&& \; \nonumber \int_X e^{-\mu\hat\varphi} |\nabla F|_{\tilde g} |\nabla {\mathcal G}_t(x_0,\cdot)|_{\tilde g} e^F \omega_X^n\le  \int_X e^{-\mu\hat\varphi} |\nabla F|_{\tilde g}(\tr_{\tilde g} g)^{\frac{1}{2}} |\nabla {\mathcal G}_t(x_0,\cdot)|_{g} e^F \omega_X^n\\
&= & \; \nonumber  \int_X e^{-\mu\hat\varphi/2} |\nabla F|_{\tilde g}Q^{\frac{1}{2}} |\nabla {\mathcal G}_t(x_0,\cdot)|_{g} e^F \omega_X^n\\
&\le & \; \nonumber C Q_{\max}^{1/2} \Big( \int_X |\nabla F|_{\tilde g}^p e^F \omega_X^n \Big)^{1/p} \cdot \Big(\int_X |\nabla {\mathcal G}_t(x_0,\cdot)|_g^{p^*} \omega_t^n \Big)^{1/p^*}\\
&\le & \; \nonumber C Q_{\max}^{1/2},
\eea
}%
where we have used $p^* = \frac{p}{p-1} < \frac{2n}{2n-1}$, so the integral on $L^{p^*}(X, \omega_t^n)$ of $|\nabla {\mathcal G}_t(x_0, \cdot)|_g$ is bounded uniformly by Theorem \ref{thm:main2}. Plugging these into \neweqref{eqn:3.20}, we get
\bea
Q_{\max}\le C Q_{\max}^{1-\eta} + CQ_{\max}^{1/2} + C, \nonumber
\eea
from which we conclude $Q_{\max}\le C$ by Young's inequality, that is $$\sup_X ( e^{-\mu\varphi} |s_D|^{2\mu \varepsilon_0} _{h_D} \tr_{\omega_X} \omega_t )\le C.$$
The proof of Theorem \ref{thm:C2} is complete by noting that $e^{-\mu\varphi}$ is uniformly bounded.

\medskip

We note that Theorem \ref{thm:C2} also applies to the  complex Monge-Amp\`ere equations with a fixed background metric.

\begin{corollary}\label{cor 3.2}
Let $\varphi$ be the solution to \neweqref{eqn:fixed MA}. Suppose $F$ satisfies the conditions \neweqref{eqn:F condition} and \neweqref{eqn:F condition 1}. Given $p>2n$, the following holds
$$\sup_X |\ddbar \varphi|_{\omega_X}^2 \le C,$$
for some constant $C>0$ depending only $n,p, \omega_X, \epsilon, N, \gamma$ and $\int_X |\nabla F|^p_{\omega_X} e^F \omega_X^n$.
\end{corollary}
The example below shows that exponent $p>2n$ in Corollary \ref{cor 3.2} is also sharp since the estimate may not hold if $p<2n$, even when other conditions are valid for $F$. The case when $p=2n$ is not clear to us.

\medskip

\noindent{\bf Example 3.2.} We take the same metrics and notations as in {\bf Example 3.1}. Fix an $a>0$ close but smaller than $1$. We can pick $1< p < 2an < 2n$, and near $0\in{\mathbb C}^n$, \neweqref{eqn:example} tells that
$$|\nabla F|^p_{\omega_{{\mathbb C}^n}} e^F \le C |z|^{-p + 2an -2n},$$
which shows that the $L^1({\mathbb{CP}}^n, \omega_{FS}^n)$-norm of $( |\nabla F|^p_{\omega_{FS}} e^F)$ is uniformly bounded (i.e. independent of $\delta>0$), while conditions \neweqref{eqn:F condition} and \neweqref{eqn:F condition 1} on $F$ hold for some $\epsilon, N, $ and $\gamma$ which are independent of $\delta$. However at the points $z$ with $|z|^2 = \delta$
$$\tr_{\omega_{{\mathbb C}^n}} ( \ddbar \varphi_\delta )\sim \delta^{a -1} \to \infty \quad\mbox{as }\delta\to 0.$$

Finally we mention an application of Corollary \ref{cor 3.2} to the regularity of complex Monge-Amp\`ere equations when $e^F$ satisfies weaker regularity than being $C^2$.  For this we need a theorem from \cite{Wa}. 
Corollary \ref{cor 3.3} below may be known to experts, but we cannot find a reference in the literature, so we include the statement and a sketched proof.
\begin{corollary}\label{cor 3.3}
Let $\varphi$ be the solution to \neweqref{eqn:fixed MA} with $F$ a Lipschitz function (i.e. $|F(x) - F(y)|\le L d_{\omega_X}(x,y)$ for some $L>0$). Then there is an $\alpha\in (0,1)$ depending on $n, \omega_X$ such that
$$\|  \varphi\|_{C^{2,\alpha}(X,\omega_X)} \le C,$$
for some constant $C>0$ depending only $n,\omega_X$, and $L$. 
\end{corollary}
\noindent{\em Proof of Corollary \ref{cor 3.3}. } By smoothing out $F$ and taking limits if necessary, we may assume $F$ is a {\em smooth} function with Lipschitz constant $\le L$. Since $\int_Xe^F\o_X^n=V$, and $F$ is Lipschitz, it follows that $\sup e^F$ and $\inf e^F>0$ are both bounded depending on $\omega_X$ and $L$, and so is the $L^p(e^F \omega_X^n)$-norm of $|\nabla F|_{\omega_X}$. It then follows from the $C^2$ estimate in Corollary \ref{cor 3.2} and the equation \neweqref{eqn:fixed MA} that $\omega_X + \ddbar \varphi$ is equivalent to $\omega_X$. We can now invoke the main theorem in \cite{Wa} to conclude the proof of Corollary \ref{cor 3.3}.

\subsection{$C^3$ estimates} To keep the notations lighter, we only consider the $C^3$ estimates for the equation \neweqref{eqn:fixed MA} with a {\em fixed} background metric. We follow closely  the approach in \cite{PSS}. 
We continue to denote $\tilde g$ and $g$ the associated metrics of $\omega_X$ and $\omega = \omega_X+\ddbar \varphi$, respectively. Furthermore we assume there is a $\theta>1$ such that
\begin{equation}\label{eqn:equivalence}
\theta^{-1} g \le \tilde g \le \theta g.
\end{equation} 
By the $C^2$ estimates in Corollary \ref{cor 3.2}, \neweqref{eqn:equivalence} holds with $\theta$ depending additionally on $\inf_X e^F$. We remark that Theorem \ref{thm:C3} below has been known with the constant $C>0$ depending on the $C^3$-norm of $F$ (cf. \cite{PSS, Y}).
\begin{theorem}\label{thm:C3}
Fix  $p>2n$. The following estimate holds:
$$\sup_X | \nabla_{\tilde g} \ddbar \varphi  |_{\tilde g}^2 \le C$$ for some $C>0$ depending on $n, p, \theta, \omega_X$, $\int_X |\nabla F|^p_{\omega_X} e^F \omega^n_X$ and $\int_X |\ddbar F|^p_{\omega_X} e^F \omega^n_X$.
\end{theorem}

\noindent{\em Proof of Theorem \ref{thm:C3}.}
As in \cite{PSS}, we denote
$$S^i_{jk} = \Gamma^i_{jk} - \tilde \Gamma^i_{jk}$$ to be
the difference of the Christoffel symbols of $g$ and$\tilde g$. Note that $S^i_{jk}$ is indeed a tensor. 
We write $$|S|_g^2 = S^i_{jk} \overline{S^l_{pq}} g_{\bar l i} g^{j\bar k} g^{k\bar q}.$$ Under the assumption \neweqref{eqn:equivalence}, it is easy to see $|S|_g^2$ is equivalent to $|\nabla_{\tilde g}\ddbar  \varphi|^2_{\tilde g}$, so it suffices to estimate $|S|_g^2$. By the calculations in \cite{PSS}, we have
\bea
\Delta_g |S|_g^2 &= &\, \label{eqn:C3} |\nabla S|_g^2 + |\bar \nabla S|_g^2 - 2 Re (\overline{S^{i}_{jk}}  R_{\bar i p, q} g^{p\bar j} g^{q \bar k})\\
&&\, \nonumber  + S*S*\ric(g) + S*S*Rm(\tilde g) + S* \nabla^{\tilde g} \ric(\tilde g),
\eea
where $T*S$ means certain linear combinations of the tensors $T$ and $S$ contracted by $g$ or $\tilde g$, and $R_{\bar i p, q}$ denotes the covariant derivative of $\ric(g)_{\bar i p}$ with respect to $g$. From the equation \neweqref{eqn:fixed MA}, we have $\ric(g) = \ric(\tilde g) - \ddbar F$. Let $G$ be the Green's function of $g$ and as before ${\mathcal G} $ be the {positive Green's function} defined in \neweqref{eqn:positive Green}. By Green's formula, for any $x\in X$, we have
%
{\small
\bea
 &&\, \nonumber  |S|_g^2(x) -\frac{1}{V}\int_X  |S|_g^2 \omega^n = \int_X {\mathcal G}(x,\cdot) \Big(  - |\nabla S|_g^2 - |\bar \nabla S|_g^2 + 2 Re (\overline{S^{i}_{jk}}  R_{\bar i p, q} g^{p\bar j} g^{q \bar k})\\
&&\qquad \nonumber  \qquad + S*S*\ric(g) + S*S*Rm(\tilde g) + S* \nabla^{\tilde g} \ric(\tilde g)   \Big)\omega^n\\
&\le &\, \nonumber  \int_X {\mathcal G}(x,\cdot) \Big(  - |\nabla S|_g^2 - |\bar \nabla S|_g^2 + C |S|_g^2 + C |S|_g^2 |\ddbar F|_{\tilde g} + C| S|_g\\
& & \, \nonumber\qquad  +  |\nabla S|_g |\ric|_g + |\bar \nabla S|_g |\ric |_g    \Big)\omega^n + \int_X 2|\nabla {\mathcal G}|_g |S|_g |\ric|_g    \omega^n \\
 &\le &\, \label{eqn:3.26}  C \int_X {\mathcal G}(x,\cdot) \Big( |S|_g^2 +  |S|_g^2 |\ddbar F|_{\tilde g} + |\ddbar F|_{\tilde g}^2 +1  \Big)\omega^n + \int_X 2|\nabla {\mathcal G}|_g |S|_g |\ric|_g    \omega^n ,
\eea }%
where in the last line we apply Cauchy-Schwarz inequality. Recall we assume a bound on $\int_X |\ddbar F|_{\tilde g}^p \omega^n$ for $p>2n$. We integrate \neweqref{eqn:YauC2} against $\omega^n$ over $X$. We obtain by H\"older inequality
\begin{equation}\label{eqn:S bound}
\int_X |S|_g^2 \omega^n \le C + \int_X |\Delta_{\tilde g} F| \omega^n\le C + \int_X |\ddbar F|_{\tilde g} \omega^n\le C.
\end{equation}
Let $x_0\in X$ be a maximum point of the function $|S|_g^2$ and denote $M = |S|_g^2(x_0)$. We apply \neweqref{eqn:3.26} at $x_0$. It follows that (denote $q=p/2>n$)
{\small
\bea
M&\le &\, \nonumber C + C M^{1-\eta}  \xk{ \int {\mathcal G}^{q^*}  }^{1/q^*} \xk{ \int |S|_g^{2\eta q}  }^{1/q} + C \xk{ \int {\mathcal G}^{q^*}   }^{1/q^*} \xk{\int |\ddbar F|_{\tilde g}^p }^{1/q} + CM^{1/2}\\
&&\, \nonumber  + C M^{1/2}\xk{ \int |\nabla {\mathcal G}|_g^{p^*}  }^{1/p^*} \xk{ \int |\ddbar F|_{\tilde g}^{p}  }^{1/p} + C M^{1-\eta} \xk{\int {\mathcal G}^{\frac{2n}{2n-1}}}^{\frac{2n-1}{2n}} \xk{\int_X |\ddbar F|_{\tilde g}^p}^{\frac{1}{p}} \xk{\int |S|_g^{2\eta b}}^{\frac{1}{b}}
\eea
}%
where $b = (p-2n)/2np$. We can choose $\eta>0$ such that $\max(2\eta q, 2\eta b) = 2$. All the integrals above are bounded, due to \neweqref{eqn:S bound} and Theorem \ref{thm:main2}. Hence we have
$$M \le C + C M^{1-\eta} + CM^{1/2},$$
from which we conclude $M\le C$. This finishes the proof of Theorem \ref{thm:C3}.



\bigskip


\noindent Department of Mathematics \& Computer Science, Rutgers University, Newark, NJ 07102

\noindent bguo@rutgers.edu

\medskip

\noindent Department of Mathematics, Columbia University, New York, NY 10027

\noindent phong@math.columbia.edu

\medskip
\noindent Department of Mathematics \& Computer Science, Rutgers University, Newark, NJ 07102

\noindent sturm@andromeda.rutgers.edu

\end{document}